\documentclass[a4paper]{article}

\usepackage[utf8]{inputenc}
\usepackage[width=15cm,height=22cm]{geometry}
\usepackage{amsmath,amsfonts,amssymb,amsthm}
\usepackage{enumitem} 
\usepackage[hidelinks,pdfusetitle]{hyperref}
\usepackage[capitalise]{cleveref}
\usepackage[dvipsnames]{pstricks}
\usepackage{bbm}
\usepackage{stmaryrd}

\crefname{equation}{}{}

\newtheorem{proposition}{Proposition}[section]
\newtheorem{theorem}[proposition]{Theorem}

\newtheorem{definition}[proposition]{Definition}
\newtheorem{corollary}[proposition]{Corollary}
\newtheorem{lemma}[proposition]{Lemma}

\newtheorem{remark}[proposition]{Remark}
\numberwithin{equation}{section}

\title{Small-time controllability of affine systems \\ and quadratic Lie brackets}

\author{Karine Beauchard\texorpdfstring{\thanks{Univ Rennes, ENS Rennes, IRMAR - UMR 6625, F-35000 Rennes, France (karine.beauchard@ens-rennes.fr)}}{},
Eugenio Pozzoli\footnote{Univ Rennes, CNRS, IRMAR - UMR 6625, F-35000 Rennes, France (eugenio.pozzoli@univ-rennes.fr)}}

\newcommand{\N}{\mathbb{N}}

\newcommand{\C}{\mathbb{C}}
\newcommand{\R}{\mathbb{R}}

\newcommand{\cE}{\mathcal{E}} 
\newcommand{\cM}{\mathcal{M}}
\newcommand{\cL}{\mathcal{L}} 

\newcommand{\ft}{\mathfrak{t}}
\newcommand{\f}{\overrightarrow}
\newcommand{\ad}{\operatorname{ad}}

\newcommand{\Lie}{\operatorname{Lie}}

\newcommand{\kb}[1]{{\color{ForestGreen}(\textbf{KB: }#1)}}

\begin{document}

\maketitle

\begin{abstract}
We prove a new sufficient condition for small-time global exact controllability of finite-dimensional control-affine systems. It points out that, under a natural Lie algebra rank condition,  quadratic Lie brackets are the only possible obstacles to controllability. This result extends to infinite-dimensional settings such as the global approximate controllability of diffeomorphisms. It also allows us to interpret recent results in PDE control theory.
\end{abstract}


\section{Introduction}

We consider an affine system on a smooth connected submanifold $\cM$ of $\R^d$
\begin{equation} \label{affine_syst}
    \dot{x}(t)=f_0(x(t))+\sum_{j=1}^{m} u_j(t) f_j(x(t))
\end{equation}
where 
$m \in \N^*$,
$f_0,\dots,f_m$ are smooth vector fields on $\cM$ and
$u=(u_1,\dots,u_m):\R_+ \to \R^m$. 

Smooth vector fields on $\cM$ are smooth vector fields on an open neighborhood of $\cM$ in $\R^d$ which are tangent to  $\cM$: for every $x \in \cM$, $f(x) \in T_x \cM$.

For every $x_0 \in \cM$ and $u \in L^1(\R_+;\R^m)$, there exists a unique solution of \eqref{affine_syst}, associated with the initial condition $x(0)=x_0$ and defined on a maximal interval, denoted  $x(t;u,x_0)$. Each time the notation $x(t;u,x_0)$ is used then, implicitly, the associated ODE is well posed on $[0,t]$.

\begin{definition}[STC]
The system \eqref{affine_syst} is small-time controllable if, for every $x_0, x_* \in \cM$, $\epsilon>0$, there exists $T \in [0,\epsilon]$ and $u \in L^1((0,T);\R^m)$ such that $x(T;u,x_0)=x_*$.
\end{definition}

The known conditions for controllability of affine systems are formulated on Lie brackets.

\subsection{Lie brackets and conditions for STC}

\begin{definition}[Lie bracket of vector fields]
For two smooth vector fields $g, h$ on $\cM$, their Lie bracket is the smooth vector field on $\cM$ defined by 
$[g,h](x)=\ad_g(h)(x):=Dh(x)g(x)-Dg(x)h(x)$
(see Lemma \ref{lem:tangent-local}).
\end{definition}

\begin{definition}[Generated Lie algebra]
For a set $F$ of smooth vector fields on $\cM$, the Lie algebra generated by $F$, denoted $\text{Lie}(F)$, is the smallest vector space of smooth vector fields on $\cM$ that contains $F$ and is stable under the Lie bracket. 
\end{definition}

The more famous sufficient condition for STC of \eqref{affine_syst} is a consequence of the Rashevski-Chow Theorem \cite{Rashevski1938, Chow1939} (see also Section \ref{subsec:Chow}).

\begin{theorem}[Rashevski-Chow] \label{thm:Chow}
Let $\cM$ be a smooth connected manifold and $F$ be a set of smooth 
vector fields on $\cM$ such that $\text{Lie}(F)(x)=T_x \cM$ for every $x \in \cM$. Then for every $x_0, x_* \in \cM$,
there exist $n \in \N^*$, $h_1,\dots,h_n \in F$, $\alpha_1, \dots,\alpha_n \in \R$ such that 
$x_*=e^{\alpha_n h_n} \dots e^{\alpha_1 h_1}(x_0)$.
\end{theorem}


The following statement is proved in \cite[page 106 for piecewise constant controls, page 110 for smooth controls]{jurdjevic-book} or  Section \ref{subsec:strong_LARC}.

\begin{theorem} \label{thm:Old_SC}
Let $\cM$ be a smooth connected manifold and 
$f_0,\dots ,f_m$ be smooth vector fields on $\cM$. 
If $\text{Lie}(\{f_1,\dots,f_m\})(x)=T_{x}\cM$ for every $x \in \cM$, then the system \eqref{affine_syst} is STC.
\end{theorem}

This sufficient condition for the STC of \eqref{affine_syst} is not necessary. A counterexample is provided by any STC system with $m<\text{dim}(\cM)$ and $[f_j,f_k]=0$ for every $j,k \in \llbracket 1 , m \rrbracket$, for instance an underactuated linear system satisfying the Kalman rank condition.
In this article, our goal is to formalize a better sufficient condition, relying on a larger Lie algebra.

\subsection{Free Lie algebra and evaluated Lie brackets}
To organize the investigation of Lie brackets in a systematic way (not depending on the system), a nice approach consists in indexing the terms by an underlying Lie bracket of the free Lie algebra. We thus introduce the following notions.

\begin{definition}[Free Lie Algebra] \label{def:FreeLie}
Let $X=\{X_0,\dots, X_m\}$ be a set of noncommutative indeterminates. For two polynomials $a,b \in \R \langle X \rangle$, we define their Lie bracket $[a,b]:=ab-ba$. For instance $[X_1,[X_1,X_0]]=X_1^2 X_0 - 2 X_1 X_0 X_1 + X_0 X_1^2$. The free Lie algebra $\cL(X)$ is the smallest vector subspace of $\R \langle X \rangle$ containing  $X$ and stable under the Lie bracket.
\end{definition}

\begin{definition}[Quadratic brackets]
For $(j,k) \in \llbracket 1,m\rrbracket \times \N^*$ we define the Lie brackets 
    \begin{equation} \label{def:MW}
        M_0^j=X_j, \qquad
        M_k^j=[M_{k-1}^j,X_0]
        \qquad \text{ and } \qquad
        W_{k}^j=[M_{k-1}^j,M_{k}^j].
    \end{equation}
For instance $W_1^j=[X_j,[X_j,X_0]]$. The bracket $W_{k}^{j}$ is called ``quadratic'' because it is of degree 2 with respect to the variable $X_j$.
\end{definition}

\begin{definition}[Evaluated Lie bracket] \label{def:evaluatedLB}
If $f_0, \dots, f_m$ are smooth vector fields on $\cM$ and 
$b\in \mathcal{L}(X)$ is an iterated Lie bracket of $X_0,\dots,X_m$,
we denote by $f_b$ the smooth vector field on $\cM$ obtained by replacing $X_j$ by $f_j$ in $b$ for  $j=0,\dots,m$. For instance, $f_{W_1^j}=[f_j,[f_j,f_0]]$. This evaluation mapping $b \mapsto f_b$ is extended linearly to $\cL(X)$ and provides a morphism of Lie algebras $\cL(X) \to \text{Lie}(f)$ where $f:=\{f_0,\dots,f_m\}$. 

If $\mathcal{E}(X)$ is a subset of $\mathcal{L}(X)$ then $\mathcal{E}(f):=\{f_b ;b \in \mathcal{E}(X) \}$ is a set of vector fields on $\cM$ and, 
for $x \in \cM$,
$\mathcal{E}(f)(x):=\{f_b(x);b \in \mathcal{E}(X) \}$ is a subset of $T_x \cM$. For instance $\cL(f)=\text{Lie}(f)$.
\end{definition}

\subsection{A natural Lie algebra rank condition for STC}

\begin{definition} \label{def:cL0(f)}
$\mathcal{L}^{\oslash}(X)$ denotes the Lie ideal of $\cL(X)$ generated by $\{X_1,\dots,X_m\}$ or, equivalently
\begin{equation} \label{def:Loslash(X)}
    \mathcal{L}^{\oslash}(X) = \text{Lie} \left(\, \left\{ \, M_k^j \, ; \, (j,k) \in \llbracket 1 , m \rrbracket \times \N \, \right\}\, \right).
\end{equation}
\end{definition}

For instance, the quadratic brackets $W_{k}^j$, defined in \eqref{def:MW}, belong to $\cL^{\oslash}(X)$. Moreover, $\cL(X)=\cL^{\oslash}(X) \oplus\R X_0$. By Definition \ref{def:evaluatedLB}, for smooth vector fields $f_0,\dots,f_m$ on $\cM$, then $\mathcal{L}^{\oslash}(f)$ is a vector subspace of $\mathcal{L}(f)$ of codimension $1$ when $f_0 \notin \mathcal{L}^{\oslash}(f)$ and $0$ otherwise. The Lie algebras $\mathcal{L}(f)$ and $\mathcal{L}^{\oslash}(f)$ play a role in characterizing the concept of accessibility.

\begin{definition}
The reachable set from $x_0 \in \cM$, at time $T \in [0,\infty)$ for the system \eqref{affine_syst} is 
$$R(x_0,T)=\left\{ \, x_* \in \cM \, ; \, \exists u \in L^1(0,T)  \text{ such that }  x(T;u,x_0)=x_* \, \right\}.$$
The reachable set from $x_0$, at time $\leq T$ is 
$R(x_0,[0,T])= \underset{T' \in [0,T]}{\cup} R(x_0,T')$.
\end{definition}

For instance, the system $\dot{x}=1$ (i.e. $\cM=\R$, $f_0=1$, $m=0$) satisfies $R(0,T)=\{T\}$ and $R(0,[0,T])=[0,T]$. 
The following equivalences are proved in \cite[Theorems 3.1 and 3.2]{Jurdjevic_Sussmann_LARC} when only piecewise constant controls are allowed (see also \cite{Hermann1963,Nagano1966, Sussmann1973}).

\begin{proposition}[Accessibility] \label{Prop:Acc}
We consider the system \eqref{affine_syst} with
$\cM$ a connected analytic manifold,
$m \in \N^*$, 
$f_0,\dots,f_m$ analytic vector fields on $\cM$ and
$x_0 \in \cM$. Then
\begin{equation} \label{LARC0}
\Big(\, \forall T>0,\, \operatorname{int} \left( R(x_0,[0,T])\, \right) \neq \emptyset
\Big) \Leftrightarrow \Big(\, \mathcal{L}(f)(x_0)=T_{x_0}\cM\, \Big),
\end{equation}
\begin{equation} \label{LARC}
\Big(\, 
\forall T>0,\, \operatorname{int} \left( R(x_0,T) \right) \neq \emptyset
\,\Big) 
\Leftrightarrow 
\Big(\, 
\cL^{\oslash} ( f )(x_0)=T_{x_0}\cM
\,\Big).
\end{equation}
\end{proposition}

The analyticity assumption is necessary for the implications "$\Rightarrow$". Conversely, the implications "$\Leftarrow$" can be proved when $\cM, f_0, \dots,f_m$ are smooth (see Section \ref{subsec:acc}).

\medskip

Setting the time $T$ for accessibility forces the vector fields $f_1,\dots,f_m$ to compensate for the drift $f_0$, just as they must do for the small-time controllability.  It is therefore natural to formulate a sufficient condition for small-time controllability based on the Lie algebra $\cL^{\oslash}(f)$. We emphasize the equivalence between $\cL^{\oslash}(f)(x)=T_{x}\cM$ and the existence of $k_*^1,\dots,k_*^m \in \N$ such that 
\begin{equation} \label{k1kn}
\begin{array}{c}
\cE(f)(x)= T_x \cM 
\qquad \text{ where } \qquad 
\cE(X)=\text{Lie} \left( \, \left\{ \, M_k^j \, ; \, (j,k) \in \llbracket 1 , m \rrbracket \times \llbracket 0 , k^j_* \rrbracket \, \right\} \, \right).
\end{array}
\end{equation}
Our main result is stated under the stronger assumption that the integers $k_*^j$ do not depend on $x$.

\subsection{STAR maps and vector fields}

To state our main result, we need the concept of small-time approximately reachable vector field.

\begin{definition}[STAR map/vector field] \label{Def:STAR}
A map $L:\text{Dom}(L) \subset \cM \to \cM$ is small-time approximately reachable if for every $x_0 \in \text{Dom}(L)$, $\epsilon>0$, there exist $T \in [0,\epsilon]$ and $u \in L^1(0,T)$ such that $|x(T;u,x_0)-L(x_0)|<\epsilon$. A vector field $g$ on $\cM$ is small-time approximately reachable if, for every $\alpha \in \R$, the flow $e^{\alpha g}$ is a STAR map.
\end{definition}

The set of STAR maps forms a semigroup, closed for the pointwise convergence (see Section \ref{subsec:semigroup}). The set of STAR vector fields has a Lie algebra structure and contains the $f_j$ so it contains also $\mathcal{L}(\{f_1,\dots,f_m\})$ (see Section \ref{subsec:LieAlg}).

\subsection{Main result}


\begin{theorem} \label{thm:Main}
Let 
$k^1_*,\dots, k^m_* \in \N \cup \{+\infty\}$,
$\cM$ be a smooth connected manifold and
$f_0, \dots,f_m$ be smooth vector fields on $\cM$.
The system \eqref{affine_syst} is STC under the following conditions:
\begin{itemize}
\item[(1)] for every $x \in \cM$, the Lie algebra rank condition \eqref{k1kn} holds,
\item[(2)] for every $(j,k) \in \llbracket 1 , m \rrbracket \times \llbracket 1 , k_*^j \rrbracket$, the vector field $f_{W_k^j}$ is STAR,
\item[(3)] for every $(j,k) \in \llbracket 1 , m \rrbracket \times \llbracket 0 , k^j_*-1 \rrbracket$, the space $\mathcal{F}_{j,k}(f)$ has finite dimension 
where
\begin{equation} \label{def:Fjk}
\mathcal{F}_{j,k}(X) = \text{Span} \{ \, \ad_{M^j_k}^{r}(X_0) \, ; \, r \in \N \, \}.
\end{equation}
\end{itemize}   
\end{theorem}

When $k_{*}^j = +\infty$, assumptions \emph{(1), (2), (3)} must hold for all $k \in \N$. This statement contains Theorem \ref{thm:Old_SC}, in the particular case $k_j^*=0$, because $M^j_0=X_j$ and the assumptions \emph{(2), (3)} are empty. But in the general case, the Lie algebra assumption \emph{(1)} involves a larger family $\mathcal{E}(X)$. Hypothesis \emph{(3)} is a technical assumption. 

\medskip

To illustrate the hypothesis \emph{(2)}, we consider the following system and associated vector fields on $\cM=\R^2$
\begin{equation} \label{Can_Syst_W1}
\begin{array}{|l}
    \dot{x}_1 = u \\
    \dot{x}_2= x_1^2 
\end{array}
\qquad \qquad
f_0(x)=\begin{pmatrix}
    0 \\ x_1^2
\end{pmatrix},
\quad
f_1(x)=\begin{pmatrix}
    1 \\ 0 
\end{pmatrix}.
\end{equation}
Clearly, this system is not STC since $x_2$ can only increase: for $T \geq 0$, $x_2(T)$ cannot reach a value $x_2^*<x_2^0$. Let us check that this system satisfies all the assumptions of  Theorem 
\ref{thm:Main} except the hypothesis \emph{(2)}:
\begin{itemize}
    \item the vector fields $f_{M_1}=[f_1,f_0]=\partial_{x_1} f_0=2x_1 e_2$ and $f_{W_1}=[f_1,f_{M_1}]=\partial_{x_1} f_{M_1} =2 e_2$ are complete and $\text{Span}\{f_1(x) , f_{W_1}(x)\}=\R^2$ thus the hypothesis \emph{(1)} is satisfied with $k_*^1=1$,
    \item for every $r \geq 3$, $\ad_{f_1}^r(f_0)=\partial_{x_1}^r f_0=0$ thus the hypothesis \emph{(3)} is also satisfied,
    \item for $\alpha \in \mathbb{R}$, we have the following explicit expressions, where $u_1(t):=\int_0^t u$
    $$\exp (\alpha f_{W_1})(x^0)=x^0 + \begin{pmatrix}
         0 \\ 2 \alpha
    \end{pmatrix},
    \qquad
    x(T;u,x^0)= x^0 + \begin{pmatrix}
        u_1(T) \\ 
        \int_0^T \left( x_1^0+u_1 \right)^2
    \end{pmatrix}
    $$
    in particular
    \begin{itemize}
    \item if $\alpha>0$, then $\exp (\alpha f_{W_1})$ is a STAR map because the above expressions are equal when $u_1(t)=\sqrt{\frac{4\alpha}{T}-2(x_1^0)^2} \sin( \frac{2\pi t}{T} )$  and $T$ is small,    
    \item if $\alpha<0$ then $\exp (\alpha f_{W_1})$ is not a STAR map because for every $T>0, x^0\in\R^2, u \in L^1$,
    $\left| x(T;u,x^0)- \exp (\alpha f_{W_1})(x^0) \right|
    \geq 2 |\alpha|$, therefore the hypothesis \emph{(2)}  is not satisfied.
    \end{itemize}
\end{itemize}
Finally, the hypothesis \emph{(2)} emphasizes that, under the Lie algebra rank condition \emph{(1)}, the quadratic brackets $W_k^j$ are the only possible obstacles to STC. 
In particular, Theorem \ref{thm:Main} implies the following statement, since
$$f_{W_k^j}=0 \qquad \Rightarrow \qquad  \mathcal{F}_{j,k}(f)=\text{Span} \left\{ \, f_0 \, , \, f_{M_{k+1}^j} \, \right\}.$$

\begin{corollary} \label{Cor:Main}
With the same notations as in Theorem \ref{thm:Main}, if (1) holds and  $f_{W_k^j}=0$ for every $(j,k) \in \llbracket 1 , m \rrbracket \times  \llbracket 1 , k^j_* \rrbracket$ then the system \eqref{affine_syst} is STC. 
\end{corollary}

\medskip

This statement proves the STC of any linear control system on $\cM=\R^d$ satisfying the Kalman rank condition. Indeed, in this case, $f_0(x)=Ax$ where $A \in \R^{d \times d}$, 
$f_j(x)=b_j \in \R^d$,
$f_{M_k^j}(x)=A^k b_j$ 
and
$f_{W_{k}^j}=0$.

\subsection{Examples}\label{sec:examples}

We illustrate Theorem \ref{thm:Main} by applying it to 3 examples. Verifying assumptions \emph{(1)} and \emph{(3)} amounts to computing iterated Lie brackets. Verifying assumption \emph{(2)} is more difficult: we have to prove the STAR property for the quadratic Lie bracket involved.  In Example 1, this quadratic bracket vanishes i.e. we apply the Corollary \ref{Cor:Main}. In Example 2,
it belongs to $\cL(\{f_1,\dots,f_m\})$  (i.e. we use Section \ref{subsec:LieAlg}). Finally, in Example 3, we prove the STAR property by hand, going back to the definition.

\medskip

\noindent \emph{Example 1:} We consider the following system and associated vector fields on $\cM=\R^3$
\begin{equation} \label{Ex1}
\begin{array}{|l}
    \dot{x}_1 = 1 + u_1 \\
    \dot{x}_2=u_1+u_2 x_1 \\
    \dot{x}_3=-x_2+u_1 x_1
\end{array}
\qquad \qquad
f_0(x)=\begin{pmatrix}
    1 \\ 0 \\ -x_2
\end{pmatrix}, \quad
f_1(x)=\begin{pmatrix}
    1 \\ 1 \\ x_1
\end{pmatrix}, \quad
f_2(x)=\begin{pmatrix}
    0 \\ x_1 \\ 0
\end{pmatrix}.
\end{equation}
A direct computation gives $[f_1,f_2](x)=e_2$. 
Thus $[[f_1,f_2],f_j](x)=\partial_{x_2} f_j(x)=0$ for $j=1,2$.
Therefore, for every $x \in \R^3$, $\cL(\{f_1,f_2\})(x)=\text{Span}\{ f_1(x), e_2\}$ i.e. Theorem \ref{thm:Old_SC} does not apply.

A direct computation gives $f_{M_1^1}(x)=[f_1,f_0](x)=-2e_3$ thus
$f_{W_1}(x)=[f_1,f_{M_1^1}](x)=2 \partial_{x_3} f_1(x) = 0$.
In particular, for every $x \in \R^3$, 
$\text{Span} \{ f_1 , [f_1,f_2] , f_{M_1^1} \}(x)=\R^3$
i.e. the assumption \emph{(1)} of Theorem \ref{thm:Main} holds with $(k_*^1,k_*^2)=(1,0)$. By Corollary \ref{Cor:Main}, the system \eqref{Ex1} is STC.

\bigskip

\noindent \emph{Example 2:} We consider the following system and associated vector fields on $\cM=\R^4$
\begin{equation} \label{Ex2}
\begin{array}{|l}
    \dot{x}_1= u_1 \\
    \dot{x}_2= x_1 \\
    \dot{x}_3= u_2 \\
    \dot{x}_4= x_1^2 + u_2 x_1 
\end{array}
\qquad \quad 
f_0(x)=\begin{pmatrix}
    0 \\ x_1 \\ 0 \\ x_1^2
\end{pmatrix}, \qquad
f_1(x)=\begin{pmatrix}
    1 \\ 0 \\ 0 \\ 0
\end{pmatrix}, \qquad
f_2(x)=\begin{pmatrix}
    0 \\ 0 \\ 1 \\ x_1
\end{pmatrix}.
\end{equation}
We have $[f_1,f_2](x)=\partial_{x_1} f_2(x) =e_4$ thus 
$[[f_1,f_2],f_j]=\partial_{x_4} f_j=0$ for $j=1,2$.
In particular, for every $x \in \R^4,$
$\mathcal{L}(\{f_1,f_2\})(x)=\text{Span}\{e_1,e_3,e_4 \}$ and Theorem \ref{thm:Old_SC} does not apply. 

We have $f_{M_1^1}=[f_1,f_0]=\partial_{x_1} f_0= e_2 + 2 x_1 e_4$, 
and $f_{W_1^1}=[f_1,f_{M_1^1}]=\partial_{x_1} f_{M_1^1}=2 e_4$ and $\ad_{f_1}^3(f_0)=[f_1,f_{W_1^1}]=0$. In particular, for every $x \in \R^4$,
$\text{Span}\{ f_1 , f_2, [f_1,f_2], f_{M_1^1} \}(x)= \R^4$ i.e.
the hypotheses \emph{(1)} and $\emph{(3)}$ hold with $(k_*^1,k_*^2)=(1,0)$. Since $f_{W_1^1} = 2 [f_1,f_2]$, it is a STAR vector field i.e. the hypothesis \emph{(2)} holds. By Theorem \ref{thm:Main}, the system \eqref{Ex2} is STC.

\bigskip

\noindent \emph{Example 3:} We consider the following system and associated vector fields on $\cM=\R^3$
\begin{equation} \label{Ex3}
\begin{array}{|l}
    \dot{x}_1 = u_1 \\
    \dot{x}_2= u_2 \\
    \dot{x}_3= x_1^2 - x_2^2
\end{array}
\qquad \quad
f_0(x)=\begin{pmatrix}
    0 \\ 0 \\ x_1^2-x_2^2
\end{pmatrix}, \quad
f_1(x)=\begin{pmatrix}
    1 \\ 0 \\ 0
\end{pmatrix}, \quad
f_2(x)=\begin{pmatrix}
    0 \\ 1 \\ 0
\end{pmatrix}.
\end{equation}
We have $[f_1,f_2]=0$ thus, for every $x \in \R^3$, $\cL(\{f_1,f_2\})(x)=\text{Span}\{e_1,e_2\}$ i.e. Theorem \ref{thm:Old_SC} does not apply. 

We have $f_{M_1^1}=[f_1,f_0]=\partial_{x_1} f_0 = 2x_1 e_3$ and $f_{W_1^1}=[f_1,f_{M_1^1}]=\partial_{x_1} f_{M_1^1} = 2e_3$ thus
$\ad_{f_1}^3(f_0)=[f_1,f_{W_1^1}]=0$. In particular, the hypotheses \emph{(1)} and $\emph{(3)}$ hold with $(k_*^1,k_*^2)=(1,0)$. 

For $\alpha>0$, $T>0$, $u_1(t)=\sqrt{\frac{2\alpha}{T}} \frac{2\pi}{T} \cos \left( \frac{2\pi}{T} t \right)$ and $u_2=0$ then $x(T;u,\cdot)=\exp (\alpha f_{W_1^1} )$. By exchanging $u_1$ and $u_2$ we obtain another control $\widetilde{u}$ such that 
$x(T;\widetilde{u},\cdot)=\exp(-\alpha f_{W_1^1})$. This proves that the vector field $f_{W_1^1}$ is STAR. By Theorem \ref{thm:Main}, the system is STC.

\subsection{Consequences and extensions}

An interesting byproduct of the proof of Theorem \ref{thm:Main} concerns the affine systems $\eqref{affine_syst}$, with state $x(t) \in \C^d$, complex-valued controls $u_j:\R \to \C$
and holomorphic vector fields $f_0, \dots, f_m: \C^d \to \C^d$.
In this case, the assumption \emph{(2)} of Theorem \ref{thm:Main} can be relaxed (see Theorem \ref{thm:C} for details). In particular, we obtain the following result, in which the base field of any vector space is $\C$.

\begin{theorem}\label{cor:C}
    If the Lie algebra $\mathcal{L}^{\oslash}(f)$ is nilpotent and $\mathcal{L}^{\oslash}(f)(x)=\C^d$ for every $x \in \C^d$ then the system \eqref{affine_syst} is STC.
\end{theorem}

As an illustration, the control system \eqref{Can_Syst_W1} with state $x(t) \in \C^2$ and control $u:\R \to \C$ is STC. The system \eqref{Ex2} with state $x(t) \in \C^4$ and control $u:\R \to \C^2$ is STC. The system \eqref{Ex3} with state $x(t) \in \C^3$ and control $u:\R \to \C^2$ is STC.

\medskip 

Theorem \ref{thm:Main} is also adapted to the control of diffeomorphisms in Section \ref{Sec:Diffeo}.

\subsection{Additional bibliographical comments}

\subsubsection{Quadratic brackets and small-time local controllability}

For scalar input systems ($m=1$), the quadratic bracket $\ad_{f_1}^2(f_0)$ is known to generate an obstruction to the small-time local controllability around an equilibrium since the 1980'. This result is due to Sussmann \cite[Proposition 6.3]{MR710995}. It is a particular case of the strongest obstruction at each even order, relying on the bracket $\ad_{f_1}^{2k}(f_0)$, proved by Stefani \cite[Theorem 1]{MR935375}, and slightly extended in \cite[Theorem 1.10]{P2}.

The obstructions based on the Lie brackets $W_k$, with $k \geq 2$, are more recent (see \cite[Theorem 2]{JDE} and \cite[Theorem 1.11]{P2}). They have also been proved for partial differential equations (see \cite{KB_ICM} and the references therein). 

It is also known that if these quadratic brackets $W_k$ are among the controlled fields, then the system is small-time locally controllable (see \cite[Theorem 1.19]{Splitting}). Theorem \ref{thm:Main} is an analogous statement for small-time global controllability.

\subsubsection{Agrachev's good brackets}


Agrachev explores issues similar to ours in the article \cite{agrachev} on the universal control system and its quotients, whose state lives in a free algebra. 

In this context, his concept of “good Lie bracket” is similar to our concept of a “STAR vector field” but it is nevertheless different because he considers only positive amplitudes $\alpha$ (for example, $X_0$ is a “good Lie bracket” in his terminology, whereas in our framework, $f_0$ may not be a STAR vector field, since the flow $e^{-f_0}$ may not be a STAR map).

Agrachev proves that the set of "good Lie brackets" for the universal control system contains all and only vector fields belonging to the closure of the set of convex combinations of $e^{-X}X_0e^{X}, X\in{\rm Lie}\{X_1,\dots,X_m\}$ (or, of course, vector fields in the closure of ${\rm Lie}\{X_1,\dots,X_m\}$). In this article, we construct a larger set of STAR vector fields, thanks to assumption \emph{(2)} -- which does not hold on the universal system, but might hold on appropriate quotients -- which is the Lie algebra generated by the ${M_k^j}, (j,k)\in \llbracket 1 , m \rrbracket \times \llbracket 0 , k^j_* \rrbracket$. 

Moreover, Proposition \ref{Prop:key} appeared also in \cite[Proposition 8]{agrachev} under stronger hypothesis (namely, the existence of an odd $\rho$ such that ${\rm ad}_g^{\rho+1}(f_0)=0$). 


\subsection{Structure of this article}

In Section \ref{sec:Main}, we prove our main result, i.e. Theorem \ref{thm:Main}. In Section \ref{sec:C}, we prove a stronger statement for complex systems and deduce Theorem \ref{cor:C}. In Section \ref{Sec:Diffeo}, we extend Theorem \ref{thm:Main} to the control of diffeomorphisms. In Section \ref{sec:PDEs} we discuss the link between Theorem \ref{thm:Main} and several small-time controllability results recently proved on PDEs. To make the article self-contained, we propose in Section \ref{Sec:App} proofs of several classical results: Proposition \ref{Prop:Acc}, Rashevski-Chow's Theorem \ref{thm:Chow}, Theorem \ref{thm:Old_SC}, estimates on flows and a conjugation formula.

\section{Proof of the main result} \label{sec:Main}
\subsection{STAR maps form a closed semigroup} \label{subsec:semigroup}

\begin{definition}
Let $\mathcal{T} \in \R_+$. A map $L:\text{Dom}(L) \subset \cM \to \cM$ is approximately reachable in time $\mathcal{T}^+$ if for every $x_0 \in \text{Dom}(L)$, $\epsilon>0$, there exist $T \in [\mathcal{T},\mathcal{T}+\epsilon]$, $u \in L^1(0,T)$ such that $|x(T;u,x_0)-L(x_0)|< \epsilon$.  
\end{definition}

In particular the STAR maps (see Definition \ref{Def:STAR}) are the maps which are approximately reachable in time $0^+$.

\begin{lemma} \label{Lem:semigroup}
\begin{enumerate}
\item For $j \in \{1,2\}$, let $\mathcal{T}_j \in \R_+$ and $L_j:\text{Dom}(L_j) \subset \cM \to \cM$ a map approximately reachable in time $\mathcal{T}_j^+$. Then $L_2 \circ L_1$ is  approximately reachable in time $(\mathcal{T}_1+\mathcal{T}_2)^+$. 
\item Let $\mathcal{T} \in \R_+$, $(\mathcal{T}_n)_{n \in \N}$ be a nonincreasing sequence converging toward $\mathcal{T}$,
$L$ and $(L_n)_{n\in\N}$ be maps on $\cM$ such that
$\text{ Dom}(L) \subset \cap_{n \in \N} \text{Dom}(L_n)$;
for every $n \in \N$, $L_n$ is approximately reachable in time $\mathcal{T}_n^+$;
and for every $x \in \text{Dom}(L)$, 
$|L_n(x)-L(x)| \to 0$. Then $L$ is approximately reachable in time $\mathcal{T}^+$. 
\end{enumerate}
\end{lemma}

Thus, the STAR maps form a semigroup which is closed for the pointwise convergence.

\begin{proof}

\emph{1.} Let $x_0 \in \text{Dom}(L_2 \circ L_1)$ and $\epsilon>0$. 
Since $L_1(x_0) \in \text{Dom}(L_2)$, there exists 
$T_2 \in [\mathcal{T}_2,\mathcal{T}_2+\epsilon/2]$ and 
$u_2 \in L^1(0,T_2)$ such that 
$|x(T_2;u_2,L_1(x_0))-L_2 \circ L_1(x_0)|<\epsilon/2$. 
By continuity of the map $x(T_2;u_2,\cdot)$ there exists $\delta>0$ such that, for every $y_0 \in \cM$ with $|y_0-L_1(x_0)|<\delta$ then $|x(T_2;u_2,y_0)-L_2 \circ L_1(x_0)|<\epsilon/2$.
Since $x_0 \in \text{Dom}(L_1)$, there exists $T_1 \in [\mathcal{T}_1,\mathcal{T}_1+\epsilon/2]$ and 
$u_1 \in L^1(0,T_1)$ such that $|x(T_1;u_1,x_0)-L_1(x_0)|<\delta$. 
Finally, $T:=T_1+T_2 \in [\mathcal{T}_1+\mathcal{T}_2,\mathcal{T}_1+\mathcal{T}_2+\epsilon]$ and 
$|x(T; u_1 \diamond u_2,x_0)- L_2 \circ L_1(x_0)| < \epsilon$.

\emph{2.} Let $x_0 \in \text{Dom}(L)$ and $\epsilon>0$. 
There exists $n \in \N$ such that 
$|L_{n}(x_0)-L(x_0)|<\epsilon/2$ and 
$\mathcal{T} \leq \mathcal{T}_{n} < \mathcal{T}+\epsilon/2$. There exists $T \in [\mathcal{T}_n,\mathcal{T}_n+\epsilon/2]$ and $u \in L^1(0,T)$ such that $|x(T;u,x_0)-L_{n}(x_0)|<\epsilon/2$. Finally 
$T \in [\mathcal{T},\mathcal{T}+\epsilon]$ and $|x(T;u,x_0)-L(x_0)|<\epsilon$.  
\end{proof}

\subsection{STAR vector fields form a Lie algebra} \label{subsec:LieAlg}


In this section, our goal is to prove the following results
(see Definition \ref{Def:STAR}). 

\begin{proposition} \label{Lem:Lie}
\begin{enumerate}
    \item For $j \in \llbracket 1, m \rrbracket$, the vector field $f_j$ is STAR.
    \item If $g,h$ are STAR vector fields then so is $g+h$.
    \item If $g,h$ are STAR vector fields then so is $[h,g]$.
\end{enumerate}
\end{proposition}

To prove Proposition \ref{Lem:Lie}, we need the following technical results.

\begin{lemma} \label{Lem_lim}
Let $g,h$ be smooth vector fields on $\cM$.  Then the following estimates hold:
\begin{itemize}
\item uniformly with respect to $x$ in any compact subset of $\text{Dom}(e^g)$
\begin{equation} \label{lim_fj}
 e^{ f_{\epsilon} + g}(x) = e^{g}(x) + \underset{\epsilon \to 0}{o}(1)
\end{equation}
under the assumption $f_{\epsilon} \to 0$ uniformly on compact subsets of $\cM$,
\item uniformy with respect to $x$ in any compact subset of 
$\text{Dom}(e^{g+h})$
\begin{equation} \label{lim_+}
\left(e^{\frac{g}{n}}e^{\frac{h}{n}} \right)^n(x) 
= e^{g+h}(x) +
\underset{n \to \infty}{O} \left( \frac{1}{n} \right),
\end{equation}
\item uniformly with respect to $x$ in any compact subset of $\text{Dom}(e^{[h,g]})$
\begin{equation} \label{lim_crochet}
    \left(e^{-\frac{g}{\sqrt{n}}} e^{-\frac{h}{\sqrt{n}}} e^{\frac{g}{\sqrt{n}}} e^{\frac{h}{\sqrt{n}}} \right)^{n}(x) = e^{[h,g]}(x) 
+ \underset{n \to \infty}{O} \left(\frac{1}{\sqrt{n}} \right).
\end{equation} 
\end{itemize}
\end{lemma}

\begin{proof}
\emph{Proof of \eqref{lim_fj}:} 
Let $K$ be a compact subset of $\text{Dom}(e^g)$.
Then $\widetilde{K}:=\left( e^{[0,1]g}(K)+\overline{B}(0,1)\right) \cap \cM$ is a compact subset of $\cM$. For $\epsilon>0$ small enough so that
$ \|f_{\epsilon} \|_{L^{\infty}(\widetilde{K})} \exp \| \nabla g \|_{L^{\infty}(\widetilde{K})}  < 1$, we deduce from Gronwall Lemma that, for every $(t,x) \in [0,1] \times K$ 
$$\left|e^{t(f_{\epsilon}+g)}(x)-e^{tg}(x)\right| \leq 
t \| f_{\epsilon} \|_{L^{\infty}(\widetilde{K})} e^{t \| \nabla g \|_{L^{\infty}(\widetilde{K})}} .$$
\\
\emph{Proof of \eqref{lim_+}:} By Lemma \ref{Lem_estim_flots}, the map
$\Phi_t(x)=e^{tg} e^{th}(x)$ satisfies
$\Phi_t(x)=x+t(g+h)(x)+O(t^2)$ as $t \to 0$,
uniformly with respect to $x$ in any compact subset of $\cM$.
We conclude by applying Lemma \ref{Lem_iteration}.
\\
\emph{Proof of \eqref{lim_crochet}:}  By Lemma \ref{Lem_estim_flots}, the map $\Phi_t(x)=e^{-\sqrt{t}g} e^{-\sqrt{t}h} e^{+\sqrt{t}g} e^{+\sqrt{t}h}$
satisfies
$\Phi_t(x)=x+t[h,g](x)+O(t^{3/2})$ as $t \to 0$,
uniformly with respect to $x$ in any compact subset of $\cM$.
We conclude by applying Lemma \ref{Lem_iteration}.
\end{proof}

\begin{proof}[Proof of Proposition \ref{Lem:Lie}]
\emph{1.} Let $\alpha \in \R$. The map $e^{\epsilon f_0 + \alpha f_j}$ is exactly reachable in time $\epsilon$. By \eqref{lim_fj}, it pointwise converges toward $e^{\alpha f_j}$, on its domain, as $\epsilon \to 0$. By Lemma \ref{Lem:semigroup}, the map $e^{\alpha f_j}$ is STAR.
\\
\emph{2.}  Let $\alpha \in \R$. By Lemma \ref{Lem:semigroup}, for every $n \in \N^*$, the map $\left(e^{\frac{\alpha}{n}g}e^{\frac{\alpha}{n}h} \right)^n$ is STAR. By \eqref{lim_+}, it pointwise converges towards $e^{\alpha(g+h)}$ on its domain. By Lemma \ref{Lem:semigroup}, $e^{\alpha(g+h)}$ is STAR.
\\
\emph{3.}   Let $\alpha \in \R$. By Lemma \ref{Lem:semigroup}, for every $n \in \N^*$, the map $\left(e^{-\frac{\alpha}{\sqrt{n}}g} e^{-\frac{1}{\sqrt{n}}h} e^{\frac{\alpha}{\sqrt{n}}g} e^{\frac{1}{\sqrt{n}}h} \right)^{n}$ is STAR. By \ref{lim_crochet}, it pointwise converges towards $e^{\alpha [h,g]}$. By Lemma \ref{Lem:semigroup}, the map  $e^{\alpha [h,g]}$ is STAR.
\end{proof}

\subsection{The key argument} 

\begin{proposition} \label{Prop:key}
Let $\rho \in \N^*$ and $g$ be a smooth vector field on $\cM$. We assume
\begin{itemize}[topsep=0cm, parsep=0cm, itemsep=0cm]
\item[(1)] $g$ and $\ad_{g}^{\rho}(f_0)$ are STAR vector fields,
\item[(2)] $\text{Span}\{ \ad_{g}^{r}(f_0) ; r \in \N\}$ has finite dimension.
\end{itemize}
Then 
\begin{itemize}[topsep=0cm, parsep=0cm, itemsep=0cm]
\item the vector field $\ad_g^{\rho-1}(f_0)$ is STAR if $(\rho-1)$ is odd,
\item the map $e^{\alpha \ad_{g}^{\rho-1}(f_0)}$ is STAR for every $\alpha>0$ when $(\rho-1)$ is even.
\end{itemize}
\end{proposition}

\begin{remark} \label{Rk:ad_f1^3=0}
If $\ad_{f_j}^3(f_0)=0$,  we can apply Proposition \ref{Prop:key} to $g=f_j$ and $\rho=3$. We obtain the following result: if $b=W_1^j$ and $\alpha >0$, then the map $e^{\alpha f_b}$ is STAR.
\end{remark}

\begin{proof}
We introduce the notation $g_r := \ad_{g}^r(f_0)$, for $r \in \N$. 
By the assumption \emph{(1)}, $g, g_{\rho}$ are STAR vector fields thus so is any linear combination of the $g_r$ for $r \geq \rho$, by Lemma \ref{Lem:Lie}. By the assumption \emph{(2)}, there exist $R \in \N$ and $a_0,\dots,a_R \in \R$ such that 
\begin{equation} \label{E0}
    E_0:=\text{Span}\{g_r;r \in \N\}=\text{Span}\{g_0,\dots,g_R\},
\end{equation}
\begin{equation} \label{CL_gr}
    g_{R+1} = a_0 g_0 + \dots +a_R g_R.
\end{equation}

\noindent \emph{Case 1: $(a_0,\dots,a_{\rho-1}) \neq 0$.} Let $r' \in \llbracket 0 , \rho-1 \rrbracket$ be the minimal integer such that $a_{r'}\neq 0$. By applying $\ad_{g}^{\rho-1-r'}$ to \eqref{CL_gr}, we obtain an expression of $g_{\rho-1}$ as a linear combination of the $g_{r}$ for $r \geq \rho$, thus $g_{\rho-1}$ is STAR. 

\medskip

\noindent \emph{Case 2: $(a_0,\dots,a_{\rho-1}) = 0$.} Then $g_{R+1}=a_{\rho} g_{\rho}+\dots+a_R g_R$. 
This relation proves that the space $\text{Span}\{g_{\rho},\dots,g_R\}$ is stable under $\ad_g$, or equivalently 
\begin{equation} \label{grho}
    E_{\rho}:=\text{Span}\{g_r;r \geq \rho\} = \text{Span}\{g_{\rho},\dots,g_R\}.
\end{equation}
From now on, $\alpha$ denotes a real number, which is positive when $(\rho-1)$ is even. For $\epsilon>0$, we define 
\begin{equation} \label{def:beta}
\beta=\beta(\epsilon):=\left( \frac{\alpha (\rho-1)!}{\epsilon}\right)^{\frac{1}{\rho-1}}.
\end{equation}
The finite-dimensional space $E_0$ (resp. $E_{\rho}$) is stable under $\ad_g$, thus $\ad_g$ is a bounded operator on $E_0$ (resp. $E_{\rho}$) and $\exp \left( \ad_{\beta g} \right)=\sum_{r=0}^{\infty} \frac{\beta^{r}}{r!} \ad_g^{r}$
is a well defined isomorphism of $E_0$ (resp. $E_{\rho}$).

\medskip

\noindent \emph{Step 1: We prove the existence of 
$h_{\epsilon} \in E_{\rho}$, $k_{\epsilon} \in E_0$ such that
\begin{equation}  \label{Step2.1}
\exp \left( \ad_{\beta g} \right)(\epsilon f_0+h_{\epsilon})=\alpha g_{\rho-1} + k_{\epsilon}
\quad \text{ and } \quad
\|k_{\epsilon}\|_{L^{\infty}_{loc}}=\underset{\epsilon \to 0}{o}(1).
\end{equation}}
By expanding the exponential and using \eqref{def:beta}, we obtain 
$$
\exp \left( \ad_{\beta g} \right)(\epsilon f_0+h_{\epsilon})
= 
\underbrace{\sum_{r=0}^{\rho-2} \frac{\beta^r}{r!} \epsilon g_r}_{=: k_{\epsilon}}
+ \alpha g_{\rho-1} + \sum_{r \geq \rho} \frac{\beta^r}{r!} \epsilon g_r + \exp \left( \ad_{\beta g} \right)(h_{\epsilon}).
$$
We define $k_{\epsilon}$ by the first term in the right hand side, then \eqref{def:beta} implies the estimate in \eqref{Step2.1}. By \eqref{grho}, the third term in the above right hand side belongs to $E_{\rho}$, therefore the sum of the third and fourth terms equals $0$ iff
$$h_{\epsilon} := \exp \left( - \ad_{\beta g} \right) \left( \sum_{r \geq \rho} \frac{\beta^r}{r!} \epsilon g_r \right) \in E_{\rho}.$$

\medskip

\noindent \emph{Step 2: We prove the map $e^{\epsilon f_0 + h_{\epsilon} }$ is approximately reachable in time $\epsilon^+$.}  This is a consequence of Lemma \ref{Lem:semigroup} and the Trotter Kato formula \eqref{lim_+} since $h_{\epsilon}$ is STAR and, for every $x \in \text{Dom}(e^{\epsilon f_0 + h_{\epsilon} })$,
$$\left( e^{\frac{\epsilon}{n} f_0 } e^{\frac{h_{\epsilon}}{n} }  \right)^{n}(x) \xrightarrow[n \to \infty]{} e^{\epsilon f_0 + h_{\epsilon} }(x).$$

\medskip

\noindent \emph{Step 3: Conclusion.}
By Lemma \ref{Lem:semigroup}, the map
$e^{-\beta g}  \, 
e^{\epsilon f_0 + h_{\epsilon}} \,
e^{\beta g}$
is approximately reachable in time $\epsilon^+$.
This map can be written $e^{F_{\epsilon}}$ where, 
$$F_{\epsilon}=\exp \left( \ad_{\beta g} \right)(\epsilon f_0+h_{\epsilon})$$
(see Lemma \ref{Lem:Conj}) and, by Step 2, $F_{\epsilon}= \alpha g_{\rho-1} + k_{\epsilon}$. By \eqref{lim_fj}, for every $x \in \text{Dom}(e^{\alpha g_{\rho-1} })$, 
$e^{F_{\epsilon}}(x) \to e^{\alpha g_{\rho-1} }(x)$ as $\epsilon \to 0$ thus, by Lemma \ref{Lem:semigroup} the map $e^{\alpha g_{\rho-1} }$ is STAR.
\end{proof}

\subsection{Proof of Theorem \ref{thm:Main}}

\noindent \emph{Step 1: For $j \in \llbracket 1 , m \rrbracket$, we prove by induction on $k \in \llbracket 0, k^j_* \rrbracket$ that the vector field $f_{M_{k}^j}$ is STAR.} For $k=0$, this holds by Corollary \ref{Lem:Lie} since $f_{M_0^j}=f_j$. Let $k \in \llbracket 0, k^j_* - 1\rrbracket$. We assume the vector field $f_{M_k^j}$ is STAR. We conclude by applying Proposition \ref{Prop:key} with $g=f_{M_k^j}$ and $\rho=2$ because
\begin{itemize}[topsep=0cm, parsep=0cm, itemsep=0cm]
\item $\ad_{g}^2(f_0)=f_{W_{k+1}^j}$ since 
$[M_k^j,[M_k^j,X_0]]=W_{k+1}^j$, 
\item $\text{Span}\{ \ad_g^r(f_0) ; r \in \N \} = \mathcal{F}_{j,k}(f)$,
\item $[g,f_0]=f_{M^j_{k+1}}$  since $[M_k^j,X_0]=M^j_{k+1}$.
\end{itemize}

\medskip

\noindent  \emph{Step 2: We prove the system is small-time approximately controllable.} Let $x_0, x_* \in \cM$. By Theorem \ref{thm:Chow}, there exists $n \in \N^*$, 
$b_1,\dots,b_n \in \{M_k^j;j \in \llbracket 1 , m \rrbracket, k \in \llbracket 0 , k^j_* \rrbracket \}$ 
and $\alpha_1, \dots, \alpha_n \in \R$
such that
$$x_* = e^{\alpha_n f_{b_n}} \dots e^{\alpha_1 f_{b_1}}(x_0).$$
By Step 1 and Lemma \ref{Lem:semigroup}, the map
$e^{\alpha_n f_{b_n}} \dots e^{\alpha_1 f_{b_1}}$ is STAR thus
$x_*$ is small-time approximately reachable from $x_0$.
\\

\noindent  \emph{Step 3: We prove that the system is small-time exactly controllable.} Let $x_0, x_* \in \cM$ and $\epsilon>0$. The assumption \emph{(1)} implies the Lie algebra rank condition $\cL(f)(x_*)=T_{x_*}\cM$ and hence ensures the accessibility of the system running backward in time $\dot{x}=-f_0(x)-\sum_{j=1}^m u_j(t) f_j(x)$. In particular, there exists an open subset $\Omega_{\epsilon}$ of $\cM$ such that, for every $x_1 \in \Omega_{\epsilon}$, there exist $T \in [0,\epsilon/2]$ and $u \in L^1(0,T)$ such that $x(T;u,x_1)=x_*$. By Step 2, there exists $T' \in [0,\epsilon/2]$, $u' \in L^1(0,T')$ such that $x_1:=x(T';u',x_0)$ belongs to $\Omega_{\epsilon}$. We get the conclusion with the time $T+T' \in [0,\epsilon]$ and the concatenated control $u \diamond u'$.

\section{Complex systems} \label{sec:C}

In this section, our goal is to prove the following statement, in which the base field of any vector space is $\C$.

\begin{theorem}  \label{thm:C}
Let 
$k^1_*,\dots, k^m_* \in \N \cup \{+\infty\}$,
$\mathcal{E}(X)$ defined as in \eqref{k1kn},
$\mathcal{F}_{j,k}(X)$ defined as in \eqref{def:Fjk} and
$f_0, \dots,f_m$ be holomorphic vector fields on $\C^d$.
The system \eqref{affine_syst} is STC under the following conditions:
\begin{itemize}
\item[(1)] for every $x \in \C^d$, $\mathcal{E}(f)(x)=\C^d$
%
\item[(2)] for every $(j,k) \in \llbracket 1 , m \rrbracket \times  \llbracket 1 , k^j_* \rrbracket$, there exists $\rho \geq 2$ such that the vector field $f_b$ is STAR for $b=\ad_{M_k^j}^{\rho}(X_0)$ and  $b=\ad_{i M_k^j}^{\rho}(X_0)$,
\item[(3)] for every $(j,k) \in \llbracket 1 , m \rrbracket \times \llbracket 0 , k^j_*-1 \rrbracket$, the space $\mathcal{F}_{j,k}(f)$ has finite dimension. 
\end{itemize} 
\end{theorem}

This statement implies Theorem \ref{cor:C} because the nilpotency of $\cL^{\oslash}(f)$ implies \emph{(2)} and \emph{(3)}. Note that the assumption \emph{(2)} is weaker in Theorem \ref{thm:C} than in Theorem \ref{thm:Main}. Theorem \ref{thm:C} is an analogue for small-time global controllability of \cite[Theorem 1.12]{Splitting} which concerns small-time local controllability. The proof of Theorem \ref{thm:C} is the same as the one of Theorem \ref{thm:Main}, except that, in Step 1, we take the conclusion for $\ell=1$ of Proposition \ref{Prop:key_C} below.

\begin{proposition} \label{Prop:key_C}
Let $\rho \in \N^*$ and $g:\C^d \to \C^d$ holomorphic. We assume
\begin{itemize}[topsep=0cm, parsep=0cm, itemsep=0cm]
\item[(1)] $g$, $ig$, $\ad_{g}^{\rho}(h_0)$, $\ad_{ig}^{\rho}(h_0)$ are STAR vector fields,
\item[(2)] $\text{Span}\{ \ad_{g}^{r}(h_0) ; r \in \N\}$ has finite dimension.
\end{itemize}
Then, for every $\ell \in \N^*$, the vector fields $\ad_{g}^{\ell}(h_0)$ and $\ad_{ig}^{\ell}(h_0)$ are STAR.
\end{proposition}

\begin{proof}
The proof is by decreasing induction on $\ell \in \llbracket 1, \rho \rrbracket$. 
The initialization for $\ell=\rho$ holds by assumption.
Let $\ell \in \llbracket 2 , \rho \rrbracket$. We assume the vector fields $\ad_{g}^{\ell}(h_0)$ and $\ad_{ig}^{\ell}(h_0)$ are STAR. If $(\rho-1)$ is odd, then $\ad_{g}^{\ell-1}(h_0)$ and $\ad_{ig}^{\ell-1}(h_0)$ are STAR by Proposition \ref{Prop:key}. If $(\rho-1)$ is even, then $i^{\rho-1}= \pm 1$ thus we just want to prove that $\ad_{g}^{\ell-1}(h_0)$ is STAR.
By Proposition \ref{Prop:key}, for every $\alpha>0$, the maps
$e^{\alpha \ad_{g}^{\rho-1}(h_0)}$ and 
$e^{\alpha \ad_{ig}^{\rho-1}(h_0)} = e^{- \alpha \ad_{g}^{\rho-1}(h_0)}$
are STAR. This gives the conclusion.
\end{proof}

\bigskip

\noindent \emph{Example 4:} We consider the following system and associated vector fields on $\cM=\C^3$
\begin{equation} \label{Ex4}
    \begin{array}{|l}
    \dot{x}_1 = 1 + v_1 \\
    \dot{x}_2=v_1+v_2 x_1 \\
    \dot{x}_3=x_2^2+ 2 v_1 x_1^2
\end{array}
\qquad \qquad
f_0(x)=\begin{pmatrix}
    1 \\ 0 \\ x_2^2
\end{pmatrix}, \quad
f_1(x)=\begin{pmatrix}
    1 \\ 1 \\ 2 x_1^2
\end{pmatrix}, \quad
f_2(x)=\begin{pmatrix}
    0 \\ x_1 \\ 0
\end{pmatrix}.
\end{equation}
Direct computations give $[f_1,f_2]=e_2$ and
$[[f_1,f_2],f_j]=\partial_{x_2} f_j=0$ for $j=1,2$
thus, for every $x \in \C^3$, $\cL(\{f_1,f_2\})(x)=\text{Span}\{ f_1(x), e_2\}$ i.e. Theorem \ref{thm:Old_SC} does not apply. 
Explicit computations give
$$f_{M_1^1}=[f_1,f_0]=(2x_2-4x_1)e_3,
\qquad
f_{W_1^1}=[ f_1 , f_{M_1^1} ]=-2 e_3,
\qquad
\ad_{f_1}^3(f_0)=[f_1,f_{W_1^1}]=2 \partial_{x_3} f_1 =0.
$$
In particular, for every $x \in \C^3$, ${\rm Span}_\C\{f_1(x),[f_1,f_2](x),f_{W_1^1}(x)\}=\C^3$ thus hypothesis \emph{(1)} holds with $(k_*^1,k_*^2)=(1,0)$. Since $\ad_{f_1}^3(f_0)=0$, the hypotheses \emph{(2)} and \emph{(3)} hold with $\rho=3$. By Theorem \ref{thm:C}, the system \eqref{Ex4} is STC.

\section{Small-time control of diffeomorphisms} \label{Sec:Diffeo}

The goal of this section is to adapt Theorem \ref{thm:Main} to the infinite-dimensional setting of the control of diffeomorphisms.
Our goal is the $L^{\infty}_{\rm loc}$-approximate controllability of the flow map  of \eqref{affine_syst} in a group of diffeomorphisms. We are interested in two groups:
\begin{itemize}
\item the group ${\rm Diff}_c^0(\cM)$ of smooth compactly supported diffeomorphisms of $\cM$ connected to the identity,
\item the group ${\rm DHam_c(\cM)}$ of Hamiltonian compactly supported diffeomorphisms of a manifold $\cM$ equipped with a symplectic form 
(in this case the vector fields $f_0,\dots,f_m$ are Hamiltonian, hence their flows are volume-preserving).
\end{itemize}
 
The $L^{\infty}_{\rm loc}$-approximate controllability in these groups implies the approximate controllability of ensembles of points (steered by the same control) because these groups act transitively on finite sets of points of $\cM$.

\subsection{$L^{\infty}_{\rm loc}$-STAR maps form a closed semigroup}

\begin{definition}[$L^{\infty}_{\rm loc}$-STAR maps]
Let $\mathcal{T} \in \R_+$. A map $L:\text{Dom}(L) \subset \cM \to \cM$ is $L^{\infty}_{\rm loc}$-approximately reachable ($L^{\infty}_{\rm loc}$-AR) in time $\mathcal{T}^+$ if for every compact subset $K$ of $\text{Dom}(L)$ and $\epsilon>0$, there exist $T \in [\mathcal{T},\mathcal{T}+\epsilon]$ and $u \in L^1(0,T)$ such that $\| x(T;u,\cdot)-L\|_{L^{\infty}(K)} < \epsilon$.  The map $L$ is $L^{\infty}_{\rm loc}$-STAR if it is $L^{\infty}_{\rm loc}$-AR in time $0^+$.
\end{definition}

Continuity is preserved by uniform convergence, so the map $L$ is continuous on $\text{Dom}(L)$

\begin{lemma} \label{Lem:semigroup_unif}
\begin{enumerate}
\item For $j \in \{1,2\}$, let $\mathcal{T}_j \in \R_+$ and $L_j:\text{Dom}(L_j) \subset \cM \to \cM$ a map $L^{\infty}_{\rm loc}$-AR in time $\mathcal{T}_j^+$. Then the map $L_2 \circ L_1$ is  $L^{\infty}_{\rm loc}$-AR in time $(\mathcal{T}_1+\mathcal{T}_2)^+$. 
\item Let $\mathcal{T} \geq 0$, $(\mathcal{T}_n)_{n \in \N}$ be a nonincreasing sequence converging toward $\mathcal{T}$,
$L$ and $(L_n)_{n\in\N}$ be maps on $\cM$ such that, 
for every $n \in \N$, $L_n$ is $L^{\infty}_{\rm loc}$-AR in time $\mathcal{T}_n^+$. We assume that any compact subset $K$ of $\text{Dom}(L)$ is contained in $\cup_{n_* \in \N} \cap_{n \geq n_*} \text{Dom}(L_n)$ and $L_n \to L$ uniformly on $K$. Then $L$ is $L^{\infty}_{\rm loc}$-AR in time $\mathcal{T}^+$. 
\end{enumerate}
\end{lemma}

\begin{proof}
\noindent \emph{1.} Let $K$ be a compact subset of $\text{Dom}(L_2 \circ L_1)$ and $\epsilon>0$. Then $K$ is a compact subset of $\text{Dom}(L_1)$, so, by continuity of $L_1$,
$\widetilde{K}:=\left( L_1(K)+\overline{B}(0,1) \right) \cap \cM$ is a compact subset of $\cM$. There exists 
$T_2 \in [\mathcal{T}_2,\mathcal{T}_2+\epsilon/2]$ and 
$u_2 \in L^1(0,T_2)$ such that 
$\| x(T_2;u_2,\cdot)-L_2 \|_{L^{\infty}(\widetilde{K})}<\epsilon/2$. In particular
\begin{equation} \label{STAR_map1}
\forall x_0 \in K, \quad
\left| x(T_2;u_2,L_1(x_0))-L_2 \circ L_1(x_0) \right|<\epsilon/2.
\end{equation} 
The map $x(T_2;u_2,\cdot)$ is continuous on the compact set $\widetilde{K}$ thus it is uniformly continuous: there exists $\delta \in (0,1)$ such that, for every $y_0, y_0' \in \widetilde{K}$ with $|y_0-y_0'|<\delta$ then $|x(T_2;u_2,y_0)-x(T_2;u_2,y_0')|<\epsilon/2$. In particular, 
\begin{equation} \label{STAR_map2}
\forall y_0 \in \cM, \, x_0 \in K
\text{ such that } |y_0-L_1(x_0)| < \delta
\text{ then }
\left|x(T_2;u_2,y_0) - x(T_2;u_2,L_1(x_0)) \right|< \epsilon/2.
\end{equation}
There exists $T_1 \in [\mathcal{T}_1,\mathcal{T}_1+\epsilon/2]$ and 
$u_1 \in L^1(0,T_1)$ such that $\| x(T_1;u_1,\cdot)-L_1 \|_{L^{\infty}(K)}<\delta$.  Finally, for $T:=T_1+T_2$ and $x_0 \in K$, by combining \eqref{STAR_map1} and \eqref{STAR_map2} with
$y_0=x(T_1,u_1,x_0)$, we obtain
$$\begin{array}{lcl}
\left| x(T; u_1 \diamond u_2,\cdot)- L_2 \circ L_1(x_0) \right|
 & \leq  &
 \left|x(T_2;u_2,x(T_1,u_1,x_0))-x(T_2;u_2,L_1(x_0)) \right| 
 \\  &  &  + \left| x(T_2;u_2,L_1(x_0)) - L_2 \circ L_1(x_0)  \right|
\leq   \epsilon.
\end{array}$$

\noindent \emph{2.} Let $K$ be a compact subset of $\text{Dom}(L)$ and $\epsilon>0$. There exists $n \in \N$ such that 
$K \subset \text{Dom}(L_n)$,
$\|L_{n}-L\|_{L^{\infty}(K)}<\epsilon/2$ and 
$\mathcal{T} \leq \mathcal{T}_{n} < \mathcal{T}+\epsilon/2$. There exists $T \in [\mathcal{T}_n,\mathcal{T}_n+\epsilon/2]$ and $u \in L^1(0,T)$ such that $\| x(T;u,\cdot)-L_{n}\|_{L^{\infty}(K)}<\epsilon/2$. Finally 
$T \in [\mathcal{T},\mathcal{T}+\epsilon]$ and $\|x(T;u,\cdot)-L\|_{L^{\infty}(K)}<\epsilon$.  
\end{proof}

\subsection{$L^{\infty}_{\rm loc}$-STAR vector fields form a Lie algebra}

\begin{definition}[$L^\infty_{\rm loc}$-STAR vector field] \
A vector field $g$ on $\cM$ is $L^\infty_{\rm loc}$-STAR if $g$ is complete and for every $\alpha \in \R$, the flow $e^{\alpha g}$ is an $L^\infty_{\rm loc}$-STAR map.
\end{definition}

\begin{proposition} \label{Prop:unif}
\begin{enumerate}
    \item For $j \in \llbracket 1, m \rrbracket$, the vector field $f_j$ is $L^{\infty}_{\rm loc}$-STAR.
    \item If $g,h$ are $L^{\infty}_{\rm loc}$-STAR vector fields then $g+h$ is $L^{\infty}_{\rm loc}$-STAR. 
    \item If $g,h$ are $L^{\infty}_{\rm loc}$-STAR vector fields then $[g,h]$ is $L^{\infty}_{\rm loc}$-STAR.
    \item Proposition \ref{Prop:key} remains valid when STAR is replaced by $L^{\infty}_{\rm loc}$-STAR. 
    \item Let $(f_n)_{n \in \N}$ be $L^{\infty}_{\rm loc}$-STAR vector fields and $f$ be a smooth vector field on $\cM$. We assume that $f_n \to f$ uniformly on compact subsets of $\cM$. Then the vector field $f$ is $L^{\infty}_{\rm loc}$-STAR.   
\end{enumerate}
\end{proposition}

\begin{proof}
In the proof of Corollary \ref{Lem:Lie},
if we replace Lemma \ref{Lem:semigroup} by Lemma \ref{Lem:semigroup_unif}, pointwise convergences by 
uniform convergences on appropriate compact sets,
guaranteed by  Lemma \ref{Lem_lim}, we obtain the first 3 statements.

For the 4th statement, we use the same proof as for Proposition \ref{Prop:key}. Indeed, by Lemma \ref{Lem_lim}, the convergences used in Steps 2 and 3 are uniform on appropriate compact sets.

Let us prove the 5th statement. Let $\alpha \in \mathbb{R}$. By Lemma \ref{Lem:semigroup_unif}, to prove that the map $e^{\alpha f}$ is $L^{\infty}_{\rm loc}$-STAR, it suffices to prove that $e^{\alpha f_n} \to e^{\alpha f}$ uniformly on compact subset of $\text{Dom}(e^{\alpha f})$, which results from \eqref{lim_fj}. 
\end{proof}

\subsection{Main result in ${\rm Diff}^0_c(\cM)$}

\begin{definition}[$L^\infty_{\rm loc}$-STAC]
System \eqref{affine_syst} is small-time $L^\infty_{\rm loc}$-approximately controllable in ${\rm Diff}^0_c(\cM)$  if for every $\phi\in {\rm Diff}^0_c(\cM)$, the map $\phi$ is $L^\infty_{\rm loc}$-STAR.
\end{definition}

\begin{theorem} \label{thm:Main_Diffeo}
Let 
$k^1_*,\dots, k^m_* \in \N \cup \{+\infty\}$,
$\mathcal{E}(X)$ defined as in \eqref{k1kn},
$\mathcal{F}_{j,k}(X)$ defined by \eqref{def:Fjk}, 
$\cM$ be a smooth connected manifold and
$f_0, \dots,f_m$ be smooth vector fields on $\cM$.
The system \eqref{affine_syst} is $L^\infty_{\rm loc}$-STAC in ${\rm Diff}^0_c(\cM)$ under the following conditions:
\begin{itemize}
\item[(1)] for every  smooth compactly supported vector field $g$ on $\cM$, there exists $(g_{n})_{n \in \N} \subset \mathcal{E}(f)$ 
such that $g_n \to g$ uniformly on compact subsets of $\cM$,
%
\item[(2)] for every $(j,k) \in \llbracket 1 , m \rrbracket \times \llbracket 1 , k_*^j \rrbracket$, the vector field $f_{W_k^j}$ is $L^{\infty}_{\rm loc}$-STAR,
\item[(3)] for every $(j,k) \in \llbracket 1 , m \rrbracket \times \llbracket 0 , k^j_*-1 \rrbracket$, the space $\mathcal{F}_{j,k}(f)$ has finite dimension. 
\end{itemize}   
\end{theorem}

The proof of Theorem \ref{thm:Main_Diffeo} relies on the following result.

\begin{proposition} \label{Prop:Thurston}
Let $\phi \in {\rm Diff}^0_c(\cM)$. There exist $p \in \N^*$ and smooth compactly supported vector fields $g_1, \dots, g_p$ on $\cM$ such that $\phi = e^{g_p} \dots e^{g_1}$.
\end{proposition}

\begin{proof}[Proof of Proposition \ref{Prop:Thurston}]
This is a consequence of the simplicity of the group ${\rm Diff}^0_c(\cM)$ proved by Thurston \cite{thurston}. Indeed, let $\psi\in {\rm Diff}^0_c(\cM)$ and $g \in {\rm Vec(\cM)}$ i.e. a smooth compactly supported vector field on $\cM$. Denote $(\psi_{\star} g)(x)=D\psi(\psi^{-1}(x))g(\psi^{-1}(x))$ the vector field obtained as the push-forward of $g$ by $\psi$ (see Definition \ref{Def:push}). Then the conjugation formula $\psi e^g\psi^{-1}=e^{\psi_{\star} g}$ shows that the subgroup $\{ e^{g_p} \dots e^{g_1}, g_1,\dots g_p\in {\rm Vec(\cM)}, p\in \N^*\}$ generated by the flows is normal in the group ${\rm Diff}^0_c(\cM)$, hence must coincide with the whole group due to the simplicity of ${\rm Diff}^0_c(\cM)$.
\end{proof}

\begin{proof}[Proof of Theorem \ref{thm:Main_Diffeo}]

\noindent \emph{Step 1: For every $(j,k) \in \llbracket 1 , m \rrbracket \times \llbracket 0, k^j_* \rrbracket$, the vector field $f_{M_{k}^j}$ is $L^{\infty}_{\rm loc}$-STAR.} The proof is the same as for Theorem \ref{thm:Main}.

\medskip
    
\noindent \emph{Step 2: Any element of $\mathcal{E}(f)$ is $L^{\infty}_{\rm loc}$-STAR.} This is a consequence of Step 1 and Lemma \ref{Lem:semigroup_unif}.

\medskip
    
\noindent \emph{Step 3: Any smooth compactly supported vector field on $\cM$ is $L^{\infty}_{\rm loc}$-STAR.} This is a consequence of the assumption \emph{(1)}, the 5th statement of Proposition \ref{Prop:unif} and Step 2.

\medskip

\noindent \emph{Step 4: Conclusion.} Let $\phi \in {\rm Diff}^0_c(\cM)$. By Proposition \ref{Prop:Thurston}, there exist $p \in \N^*$ and smooth compactly supported vector fields $g_1, \dots, g_p$ on $\cM$ such that $\phi = e^{g_p} \dots e^{g_1}$. By Step 2, the maps $e^{g_k}$ are $L^{\infty}_{\rm loc}$-STAR. By Lemma \ref{Lem:semigroup_unif}, the map $\phi$ is  $L^{\infty}_{\rm loc}$-STAR.
\end{proof}

Some bibliographical comments are in order: related versions of Theorem \ref{thm:Main_Diffeo} corresponding to the driftless case (i.e. $f_0=0$ and $k_*^j=0$) already appeared in \cite[Theorem D]{heintze-liu} in the more general setting of Banach infinite-dimensional manifolds, and were specialized to diffeomorphism groups in \cite{agrachev-sarychev2020,agrachev-sarychev}.

\subsection{Main result in ${\rm DHam_c(\cM)}$}

In this section we suppose that the manifold $\cM$ is equipped with a symplectic form $\omega$. The notion of Hamiltonian vector field $\f h $ on $\cM$ associated to a smooth function $h\in \mathcal{C}^\infty(\cM,\R)$ is defined by the formula $\omega(\f h, V)=V(h)$ for every vector field $V$ on $\cM$. $\mathcal{C}^\infty(\cM,\R)$ then becomes a Lie algebra with bracket $\{h_1,h_2\}:=\omega(\f {h_1},\f {h_2})$, usually called Poisson bracket. Since $[\f {h_1},\f {h_2}]=\f{\{h_1,h_2\}}$, the set of Hamiltonian vector fields on $\cM$ forms a Lie subalgebra of the smooth vector fields on $\cM$. 

We moreover suppose that $f_0,\dots,f_m$ are Hamiltonian vector fields on $\cM$, hence there exist $h_0,\dots,h_m\in \mathcal{C}^\infty(\cM)$ such that $\f{h_j}=f_j, j=0,\dots,m$.

Let $h=\{{h}_0,\dots,{h}_m\}\subset \mathcal{C}^\infty(\cM)$, and we denote by ${\rm Poisson}(h)\subset \mathcal{C}^\infty(\cM)$ the Lie algebra generated by $h$ w.r.t. the Poisson bracket,  and by $\mathcal{L}(\f h)$ the Lie algebra generated by the Hamiltonian vector fields $\f{h_0},\dots,\f{h_m}$. In analogy to the notions introduced in Definition \ref{def:evaluatedLB} on evaluated brackets, if $h_0, \dots, h_m$ are smooth functions on $\cM$ and 
$b\in \mathcal{L}(X)$ is an iterated Lie bracket of $X_0,\dots,X_m$,
we denote by $h_b$ the smooth function on $\cM$ obtained by replacing $X_j$ by $h_j$ in $b$ for  $j=0,\dots,m$. For instance, $h_{W_1^j}=\{h_j,\{h_j,h_0\}\}$. This evaluation mapping $b \mapsto h_b$ is extended linearly to $\cL(X)$ and provides a morphism of Lie algebras $\cL(X) \to \text{Poisson}(h)$. 

If $\mathcal{E}(X)$ is a subset of $\mathcal{L}(X)$ then $\mathcal{E}(h):=\{h_b ;b \in \mathcal{E}(X) \}$ is a set of smooth functions on $\cM$.

\begin{definition}[$L^\infty_{\rm loc}$-STAC]
System \eqref{affine_syst} is small-time $L^\infty_{\rm loc}$-approximately controllable in ${\rm DHam}_c(\cM)$ if for every $\phi\in {\rm DHam}_c(\cM)$, the map $\phi$ is $L^\infty_{\rm loc}$-STAR.
\end{definition}

\begin{theorem} \label{thm:Main_Diffeo_Ham}
Let 
$k^1_*,\dots, k^m_* \in \N \cup \{+\infty\}$,
$\mathcal{E}(X)$ defined as in \eqref{k1kn},
$\mathcal{F}_{j,k}(X)$ defined as in \eqref{def:Fjk}, 
$\cM$ be a smooth connected manifold equipped with a symplectic form, and
$h_0, \dots,h_m$ be smooth functions on $\cM$.
The system \eqref{affine_syst} is $L^\infty_{\rm loc}$-STAC in ${\rm DHam}_c(\cM)$ under the following conditions:
\begin{itemize}
\item[(1)] for every  smooth compactly supported Hamiltonian vector field $\f s$ on $\cM$, there exists $(s_{n})_{n \in \N} \subset \mathcal{E}(h)$ 
such that $\f{s_n} \to \f s$ uniformly on compact subsets of $\cM$,
%
\item[(2)] for every $(j,k) \in \llbracket 1 , m \rrbracket \times \llbracket 1 , k_*^j \rrbracket$, the Hamiltonian vector field $\f{h_{W_k^j}}$ is $L^{\infty}_{\rm loc}$-STAR,
\item[(3)] for every $(j,k) \in \llbracket 1 , m \rrbracket \times \llbracket 0 , k^j_*-1 \rrbracket$, the space $\mathcal{F}_{j,k}(h)$ has finite dimension.
\end{itemize}   
\end{theorem}
The proof of Theorem \ref{thm:Main_Diffeo_Ham} is  analogous to the one of Theorem \ref{thm:Main_Diffeo}, the only difference being the replacement of Proposition \ref{Prop:Thurston} by the following statement.

\begin{proposition}  \label{Banyaga}
Let $\phi \in {\rm DHam}_c(\cM)$. There exist $p \in \N^*$ and smooth compactly supported functions $h_1, \dots, h_p$ on $\cM$ such that $\phi = e^{\f {h_p}} \dots e^{\f {h_1}}$, where $\f {h_i}$ denotes the Hamiltonian vector field on $\cM$ associated to $h_i$.
\end{proposition}

\begin{proof}
This is a consequence of the simplicity of the group ${\rm DHam}_c(\cM)$ proved by Banyaga \cite{banyaga}. Indeed, if $g=\f h$ for some $h\in C^\infty_c(\cM,\R)$ and $\det D \psi=1$, the conjugation formula $\psi e^{\f h}\psi^{-1}=e^{\f {h\circ\psi}}$ shows that the subgroup $\{ e^{\f {h_p}} \dots e^{\f {h_1}}, h_1,\dots h_p\in C^\infty_c(\cM), p\in \N^*\}$ generated by the flows of Hamiltonian vector fields is normal in the group ${\rm DHam}_c(\cM)$, hence must coincide with the whole group due to the simplicity of ${\rm DHam}_c(\cM)$.
\end{proof}

Some bibliographical comments are in order: approximate controllability in ${\rm DHam}(T^*\R^d)$ (where the manifold $\cM=T^* \R^d$ with global coordinates $(q_1,\dots,q_d,p_1,\dots,p_d)$ is equipped with the canonical symplectic form $\omega=\sum_{i=1}^ddq_i\wedge dp_i$) by flows along Hamiltonian vector fields was also studied in \cite{turaev} and \cite{kazandjian-pozzoli-sigalotti}, for typical examples of Hamiltonian systems in the form of kinetic plus potential energy. Approximate decomposition results for Hamiltonian diffeomorphisms (in the driftless case $h_0=0$) of $T^* \mathbb{T}^d$ and $T^*\R^d$ are also obtained in \cite{berger} .

\subsection{Example}

In this section we give an example satisfying the conditions of Theorem \ref{thm:Main_Diffeo_Ham}. Let $\cM=T^*\R^d$, with global position-momentum coordinates $(q_1,\dots,q_d,p_1,\dots,p_d)=(q,p)\in \R^d_q\times \R^d_p$, equipped with the canonical symplectic form $\sum_{i=1}^ddq_i\wedge dp_i$, inducing the canonical Poisson bracket $\{h,s\}=\sum_{i=1}^d\partial_{q_i}h\partial_{p_i} s-\partial_{q_i}s\partial_{p_i} h$ for every $h,s\in \mathcal{C}^\infty(T^*\R^d,\R), h=h(q,p),s=s(q,p)$. Given $h\in \mathcal{C}^\infty(T^*\R^d,\R)$, its corresponding Hamiltonian vector field is thus $\f{h}=\{h,\cdot\}=\sum_{i=1}^d\partial_{q_i}h\partial_{p_i} -\partial_{p_i} h \partial_{q_i}$. Let $h_0(q,p)=|p|^2+V(q)$, $h_1(q,p)=q_1,\dots,h_d(q,p)=q_d, h_{d+1}(q,p)=e^{-|q|^2/2},h_{d+2}(q,p)=e^{-|p|^2/2}$. We have 
\begin{equation}\label{eq:poisson-brackets}
\{h_0,h_j\}=p_j,\quad \{\{h_0,h_j\},h_j\}=1,\quad {\rm ad}_{h_j}^3(h_0)=0,\quad j=1,\dots,d. 
\end{equation}
Since $\{\{h_0,h_j\},h_j\}$ is constant, its associated Hamiltonian vector field $\f{{h}_{W^j_{1}}}$ is equal to $0$, hence $\f{{h}_{W^j_{1}}}$ is trivially $L^\infty_{\rm loc}$-STAR for all $j=1,\dots,d$. Moreover $\mathcal{F}_{j,0}(h)$ is two-dimensional (by \eqref{eq:poisson-brackets}) for all $j=1,\dots,d$. We set $k_*^1=\dots=k_*^d=1$ and $k_*^{d+1}=k_*^{d+2}=0$. Notice that, for $n_1,\dots,n_d\in\N$,
$${\rm ad}_{h_d}^{n_d}\dots{\rm ad}_{h_1}^{n_1}(h_{d+2})=\partial_{p_d}^{n_d}\dots\partial_{p_1}^{n_1}e^{-|p|^2/2},\quad {\rm ad}_{h_{M_1^d}}^{n_d}\dots{\rm ad}_{h_{M_1^1}}^{n_1}(h_{d+1})=\partial_{q_d}^{n_d}\dots\partial_{q_1}^{n_1}e^{-|q|^2/2}.$$
Those latter brackets span the Hermite functions in the $p$-variables (and in the $q$-variables), hence every element of $\mathcal{C}^\infty_c(\R^d_p,\R)$ can be $W^{1,\infty}$-approximated by a sequence in $${\rm span}\{{\rm ad}_{h_d}^{n_d}\dots{\rm ad}_{h_1}^{n_1}(h_{d+2}),n_1,\dots,n_d\in\N\}$$ (and every element of $\mathcal{C}^\infty_c(\R^d_q,\R)$ can be $W^{1,\infty}$-approximated by a sequence in $${\rm span}\{{\rm ad}_{h_{M_1^d}}^{n_d}\dots{\rm ad}_{h_{M_1^1}}^{n_1}(h_{d+1}),n_1,\dots,n_d\in\N\}.)$$ As a consequence, every $\f{h},h\in \mathcal{C}^\infty_c(\R^d_q,\R),$ and $\f{s},s\in \mathcal{C}^\infty_c(\R^d_p,\R)$, is $L^\infty_{\rm loc}$-STAR. Finally, for every $h\in \mathcal{C}^\infty_c(\R^d_q,\R),s\in \mathcal{C}^\infty_c(\R^d_p,\R)$, the formula $\{h,s\}=\sum_{j=1}^d\partial_{q_j}h\partial_{p_j}s$ implies that every Hermite function in the $(q,p)$-variables belongs to $\mathcal{E}(h)$. Thus, every element of $\mathcal{C}^\infty_c(T^*\R^d,\R)$ can be $W^{1,\infty}$-approximated by a sequence in $\mathcal{E}(h)$, and property (1) of Theorem \ref{thm:Main_Diffeo_Ham} follows.

\section{Discussion about PDEs} \label{sec:PDEs}

In this section, we point out connections between Theorem \ref{thm:Main} and the small-time global approximate controllability results recently obtained in \cite{duca-nersesyan,pozzoli-wave,beauchard-pozzoli1,beauchard-pozzoli2,beauchard-perrin-pozzoli} for bilinear PDEs such as Schrödinger and wave equations.

\subsection{Schr\"odinger equation}

Consider the PDE
\begin{equation}\label{eq:schro}
i\partial_t \psi(x,t)=\left(-\Delta+V_0(x)+\sum_{j=1}^mu_j(t)V_j(x)\right)\psi(x,t),\quad (x,t)\in M\times [0,T],
\end{equation}
where $M=\mathbb{T}^d$ or $\R^d$ and
$V_0, V_1,\dots,V_m\in \mathcal{C}^\infty(M,\R)$.
The controls are $u_1,\dots,u_m:\R \to \R$ and the state $\psi(t,.)$ lives in the $L^2(M;\C)$-sphere. The small-time global $L^2$-approximate controllability of such systems is proved in \cite{beauchard-pozzoli2} in 2 configurations
corresponding to $M=\mathbb{T}^d$ or $\R^d$ and appropriate choices of potentials $V_j$.

\medskip

In this proof, the first ingredient is the small-time controllability of phases: for every $\varphi \in L^2(M;\R)$ the map $\psi_0 \mapsto e^{i \varphi} \psi_0$ is small-time $L^2$-approximately reachable ($L^2$-STAR), see \cite[Proposition 22]{beauchard-pozzoli2}. The proof of this result consists of repeatedly applying the following argument 
$$\big( \forall \alpha \in \R,\, e^{i \alpha \varphi} \text{ is } L^2-\text{STAR} \big) \quad  \Rightarrow  \quad  \big(\forall \alpha \in (0,\infty), \, e^{-i\alpha |\nabla \varphi|^2} \text{ is } L^2-\text{STAR} \big).$$
This implication corresponds to Remark \ref{Rk:ad_f1^3=0}. Indeed, heuristically, with $f_0=i(\Delta-V_0)$ and $f_1=i\varphi$ 
\begin{equation} \label{crochets_Schro}
\ad_{f_1}(f_0) = - 2 \nabla\varphi \cdot \nabla -\Delta \varphi,
\qquad \qquad
\ad_{f_1}^2(f_0)= - 2i|\nabla \varphi|^2,
\qquad \qquad
\ad_{f_1}^3(f_0)=0.
\end{equation}

\medskip

The second ingredient is the small-time $L^2$-approximate reachability of transport operators associated with gradient vector fields: for every $\varphi \in \mathcal{C}^{\infty}(M;\R)$ the map $\psi_0 \mapsto e^{ \mathcal{T}} \psi_0$ is $L^2$-STAR, where $\mathcal{T}=2 \nabla \varphi \cdot \nabla + \Delta \varphi$ (see \cite[Theorem 16]{beauchard-pozzoli2}). This corresponds to Proposition \ref{Prop:key} with $g=i \varphi$ and $\rho=2$. Indeed, by the small-time controllability of phases, the maps $e^{i\alpha g}$ and $e^{i\alpha|\nabla \varphi|^2}$ are $L^2$-STAR for every $\alpha \in \R$, i.e. with the notations of \eqref{crochets_Schro} the``vector fields'' $f_1$ and $\ad_{f_1}^2(f_0)$ are STAR, thus so is the ``vector field'' $\ad_{f_1}(f_0)$ which is a transport operator.

\subsection{Wave equation}

Consider the Klein-Gordon wave equation  
\begin{equation}\label{eq:wave_bis}
{\partial^2_t}w(t,x)=\left(\Delta-1+\sum_{j=1}^{m} u_j(t)V_j(x)\right)w(t,x),\quad (t,x)\in \mathbb{R} \times \mathbb{T}^d,
\end{equation}
where 
$d \in \mathbb{N}^*$ and
$V_1,\dots,V_m \in \mathcal{C}^{\infty}(\mathbb{T}^d,\R)$.
The controls are
$u_0,\dots,u_{m}:\mathbb{R} \to \mathbb{R}$ and the state is $(w,\partial_t w )(t)$.
We introduce its vectorised version, for $W=(w,\partial_t w)$,
\begin{equation} \label{wave_vect}
\frac{dW}{dt}(t)=\left( A+\sum_{j=1}^{m} u_j(t) V_j B \right) W(t)
\quad \text{ where } \quad 
A=\begin{pmatrix}
    0 & 1 \\ \Delta-1 & 0
\end{pmatrix}
\quad \text{ and } \quad
B=\begin{pmatrix}
    0 & 0 \\ 1 & 0
\end{pmatrix}.
\end{equation}
In \cite{beauchard-perrin-pozzoli}, with appropriate potentials $V_1,\dots,V_m$, we prove several $H^1 \times L^2$-approximate controllability results for this system.

In the proof, the first ingredient is the small-time controllability of the velocity: for every $\varphi \in \mathcal{C}^{\infty}(\mathbb{T}^d)$, the map $\exp(\varphi B):(w_0,\dot{w}_0) \mapsto (w_0,\dot{w}_0+\varphi w_0)$ is small-time $H^1 \times L^2$-approximately reachable (STAR); see \cite[Proposition 24]{beauchard-perrin-pozzoli}. The proof of this result consists in iterating the following argument
$$\big( \forall \alpha \in \R,\, \exp( \alpha  \varphi B) \text{ is } H^1 \times L^2-\text{STAR} \big) \quad  \Rightarrow  \quad  \big(\forall \alpha \in (0,\infty), \, \exp( - \alpha \varphi^2 B ) \text{ is } H^1 \times L^2-\text{STAR} \big).$$
This implication corresponds to Remark \ref{Rk:ad_f1^3=0}. Indeed, heuristically, 
$$M:=[B,A]=\begin{pmatrix}
     -1 & 0 \\ 0 & 1
 \end{pmatrix}, \qquad
[\varphi B,[\varphi B,A]]=-2\varphi^2B, \qquad
\text{ad}_{B}^j(A)=0\,, \forall j\geq 3.$$

Another ingredient is the small-time $H^1 \times L^2$-approximate reachability of the map $\exp(\alpha M): (w_0,\dot{w}_0) \mapsto (e^{-\alpha} w_0, e^{\alpha }\dot{w}_0)$, for every $\alpha \in \R$ (see \cite[Proposition 27]{beauchard-perrin-pozzoli}). This corresponds to Proposition \ref{Prop:key} with $\rho=2$. By the small-time controllability of the velocity, the maps $\exp(\alpha B)$ and $\exp(-2 \alpha^2 B)$ are STAR thus so are the $\exp(\alpha M)$.

\appendix
\section{Appendix}\label{Sec:App}

In this section, we give self-contained proofs of several classical results used in this article: 
\begin{itemize}[topsep=0cm, parsep=0cm, itemsep=0cm]
\item in Section \ref{subsec:Tgt}, we prove the Lie algebra structure of the set of smooth vectorfields on $\cM$,
\item in Section \ref{subsec:inv_arg}, a local inversion argument, which is the source of all the following results
\item in Section \ref{subsec:acc}, accessibility, for smooth vector fields, under the Lie algebra rank condition,
\item in Section \ref{subsec:Chow}, Rashevski-Chow's Theorem \ref{thm:Chow},
\item in Section \ref{subsec:strong_LARC}, Theorem \ref{thm:Old_SC},
\item in Section \ref{subsec:estim_flow}, estimates on flows,
\item in Section \ref{subsec:conj}, a conjugation formula.
\end{itemize}

\subsection{Lie algebra structure} \label{subsec:Tgt}

\begin{lemma}
    \label{lem:tangent-local}
    Let $\cM$ be a smooth submanifold of $\R^d$,
    $V$ be an open subset of $\R^d$, 
    and $f, g: V \to \R^d$ be smooth vector fields tangent to $\cM$ i.e. for every $x \in \cM \cap V$, $f(x),g(x)\in T_x \cM$.
    Then $[f,g]$ is tangent to $\cM$.
\end{lemma}

\begin{proof}
    Fix $x \in \cM$. 
    Since $\cM$ is an embedded submanifold, there exist an open neighborhood $W \subset V$ of $x$ and a smooth submersion $\psi : W\to \R^{d-k}$, where $k=\dim \cM$, such that $\cM \cap W = \psi^{-1}(0)$.
    Moreover, for all $y \in \cM \cap W$, $T_y \cM = \ker D\psi(y)$.
	
	Since $g(x)\in T_x \cM$, there exists a smooth curve $\gamma:(-\varepsilon,\varepsilon) \to \cM \cap W$ such that
	$\gamma(0)=x$ and $\gamma'(0)=g(x)$. 
	Since $f$ is tangent to $\cM$, $D\psi(\gamma(t)) f(\gamma(t)) = 0$ for all $t$.
	Differentiating with respect to $t$ at $t=0$ yields
	\begin{equation}
		D^2\psi(x)\bigl(g(x),f(x)\bigr)+D\psi(x)Df(x)g(x)=0.
	\end{equation}
	Swapping $f$ and $g$, and using the symmetry of $D^2\psi(x)$, we get
	\begin{equation}
	D\psi(x) [f,g](x) = D\psi(x)\bigl(Dg(x)f(x)-Df(x)g(x)\bigr)=0.
	\end{equation}
    Thus $[f,g](x) \in T_x \cM$.
\end{proof}

\subsection{An inversion argument} \label{subsec:inv_arg}

 The goal of this section is to present the inversion argument \cref{p:access} which is a classical result of control theory, constructing a ``normal control'' as in \cite[page 169]{MR872457}. 
\begin{proposition}
    \label{p:access}
    Let $V \subset \R^d$ be an open neighborhood of $x_0 \in \R^d$, and $F$ a set of smooth vector fields on $V$.
    Assume that $\text{Lie}(F)(x_0) = \R^d$.
    Then, for all $T>0$, there exist $h_1,\dots,h_d \in F$,
    and an open subset $\Omega \subset (0,T)^d$ such that the map
    \begin{equation}
        \Phi : \ft = (t_1, \dotsc, t_d) \mapsto e^{t_d h_d} \dotsb e^{t_1 h_1}(x_0)
    \end{equation}
    is a smooth embedding of $\Omega$ into $V$.
    In particular, $D \Phi$ has rank $d$ on $\Omega$.
\end{proposition}

\begin{proof}
    Choose $g_1,\dots,g_d\in \text{Lie}(F)$ such that $g_1(x_0),\dots,g_d(x_0)$ is a basis of $\R^d$.
    By continuity, there exists an open neighborhood $W \subset V$ of $x_0$ such that, for every $x \in W$, the vectors $g_1(x),\dots,g_d(x)$ are still linearly independent, hence $\text{Lie}(F)(x) = \R^d$.

    We prove by induction on $k\in\{1,\dots,d\}$ the following statement:
    \emph{There exist $h_1,\dots,h_k\in F$ and an open subset $\Omega_k \subset (0,T)^k$ such that the map $\Phi_k : \ft = (t_1, \dotsc, t_k) \mapsto e^{t_k h_k} \dotsb e^{t_1 h_1}(x_0)$ is a smooth embedding of $\Omega_k$ into $W$.}

    \medskip
    \noindent\emph{Initialization: $k=1$.}
    There exists $h_1 \in F$ such that $h_1(x_0)\neq 0$,
    otherwise, all elements of $\Lie \mathcal{F}$ would vanish at $x_0$.
    Then $\Phi_1$ is well-defined and smooth on a neighborhood of $0 \in \R$. Since $D \Phi_1(0) = f_1(x_0) \neq 0$, by continuity, there exists an open subset $\Omega_1 \subset (0,T)$ on which $\Phi_1$ is an embedding.

    \medskip
    \noindent\emph{Induction step.}
    Assume the statement proved for some $k\in\{1,\dots,d-1\}$.

    Since $\Phi_k$ is a smooth embedding, $\cM_k := \Phi_k(\Omega_k)$ is a smooth embedded submanifold of $\R^d$ of dimension $k$.
    Moreover, for any $\mathfrak{s} \in \Omega_k$ and $y = \Phi_k(\mathfrak{s}) \in \cM_k$, $T_y \cM_k = D \Phi_k(\mathfrak{s}) (\R^k)$.

    We use the method of proof by contradiction.    
    Suppose that any field in $F$ is tangent to $\cM_k$.
    Then, by \cref{lem:tangent-local}, every field in $\text{Lie}(F)$ is tangent to $\cM_k$.
    Therefore, for every $x \in \cM_k$, $\text{Lie}(F)(x) \subset T_x \cM_k$, which is impossible since $\text{dim}(\text{Lie}(F)(x))=d>k=\text{dim}(T_x \cM_k)$.
    Hence there exist $h_{k+1} \in F$ and $y \in \cM_k$ such that $h_{k+1}(y) \notin T_y \cM_k$.
    
    Write $y = \Phi_k(\mathfrak{s})$ with $\mathfrak{s} \in \Omega_k$.
    For $\ft = (t_1, \dots, t_{k+1})$ in a neighborhood of $(\mathfrak{s},0)$, the map $\Phi_{k+1}$ is well-defined, smooth, with values in $W$.
    Moreover,
    \begin{equation}
        D\Phi_{k+1}(\mathfrak{s},0)(h,\lambda)
        =
        D\Phi_k(\mathfrak{s})h+\lambda h_{k+1}(y).
    \end{equation}
    Since $h_{k+1}(y) \notin D\Phi_k(\mathfrak{s})(\R^k)$, the map
    $D\Phi_{k+1}(\mathfrak{s},0)$ has rank $k+1$.
    Thus there exists an open subset $\Omega_{k+1} \subset \Omega_k \times (0,T)$ on which $\Phi_{k+1}$ is a smooth embedding.
\end{proof}

\subsection{Accessibility} \label{subsec:acc}

\begin{proposition} \label{Prop:Access_0}
Let $\cM$ be a smooth submanifold of $\R^d$,
$f_0, \dots, f_m$ be smooth vector fields on $\cM$ and
$x_0 \in \cM$.
We assume $\cL(f)(x_0)=T_{x_0}\cM$.
Then, for every $T>0$, $\text{Int} \big( R(x_0,[0,T]) \big) \neq \emptyset$.
\end{proposition}

\begin{proof}
We are going to prove a local result, so we can assume that $\cM=\R^d$.
We define $F:=\{f_0 \pm f_j;j \in \llbracket 1 , m \rrbracket\}$. Then $\text{Lie}(F)(x_0)=\R^d$.
By \cref{p:access}, there exists $h_1,\dots,h_d \in F$ and
an open set $\Omega \subset (0,\frac{T}{d})^d$ such that the map
 \begin{equation}
        \Phi : \ft = (t_1, \dotsc, t_d) \mapsto e^{t_d h_d} \dotsb e^{t_1 h_1}(x_0)
    \end{equation}
    is a smooth embedding of $\Omega$. In particular, the open set $\Phi(\Omega)$ is contained in $R(x_0,[0,T])$.
\end{proof}

\begin{proposition} \label{Prop:Access_1}
Let $\cM$ be a smooth submanifold of $\R^d$,
$f_0, \dots, f_m$ be smooth vector fields on $\cM$ and
$x_0 \in \cM$.
We assume $\cL^{\oslash}(f)(x_0)=T_{x_0}\cM$.
Then, for every $T>0$, $\text{Int} \big( R(x_0,T) \big) \neq \emptyset$.
\end{proposition}

\begin{proof}
We are going to prove a local result, so we can assume that $\cM=\R^d$.
We introduce the vector fields $g_0, g_1, \dots, g_m$ on $\R^{d+1}$ defined, for $y=(t,x) \in \R \times \R^d$ by
$$g_0(y)=\begin{pmatrix}
    1 \\ f_0(x)
\end{pmatrix},
\qquad 
g_j(y)=\begin{pmatrix}
    0 \\ f_j(x)
\end{pmatrix}
\text{ for } j \in \{1,\dots,m\}.$$
Let $y_0=(0,x_0) \in \R^{d+1}$.
Then the solution $y(t;u,y_0)$ of the system
\begin{equation} \label{syst_avec_t'=1}
\dot{y}=g_0(y)+\sum_{j=1}^m u_j(t) g_j(y)
\end{equation}
associated with the initial condition $y(0)=y_0$ is related to the solution $x(t;u,x_0)$ of \eqref{affine_syst} associated with the initial condition $x(0)=x_0$ by 
$y(t;u,y_0)=(t,x(t;u,x_0))$.
Explicit computations prove that, for any Lie bracket $b \in \cL(X) \setminus \{X_0\}$ then
$F_b(y)=(0 , f_b(x))$. Taking into account the explicit expression of $g_0$, the assumption $\cL^{\oslash}(f)(x_0)=\R^d$ implies
$\mathcal{L}(g)(y_0)=\R^{d+1}$. 

For $T\geq 0$, let $R(x_0,T)$ be the reachable set in time $T$ for system \eqref{affine_syst} and $R(y_0,[0,T])$ be the reachable set in time $\leq T$ for system \eqref{syst_avec_t'=1}.  

By \cref{Prop:Access_0}, one may consider an interior point $y_1=(T_1,x_1)$ of $R(y_0,[0,T])$. Then $T_1 \in [0,T]$ and there exists $r>0$ such that $\{T_1\} \times B(x_1 ,r) \subset R(y_0,[0,T])$. Then $B(x_1 ,r) \subset R(T_1,x_0)$. Since $e^{(T-T_1)f_0}$ is a diffeomorphism, then $e^{(T-T_1)f_0}B(x_1 ,r)$ is an open set of $\R^d$ contained in $R(x_0,T)$.  
\end{proof}

\subsection{Rashevski-Chow's theorem} \label{subsec:Chow}

\begin{proof}[Proof of Theorem \ref{thm:Chow}]
For $x_0 \in \cM$, we define
$$\mathcal{O}(x_0):=\left\{ e^{\alpha_n h_n} \dots e^{\alpha_1 h_1}(x_0);
n \in \N, h_1, \dots, h_n \in F, \alpha_1, \dots, \alpha_n \in \R \right\}.$$
The vector fields of $F$ are not assumed to be complete, thus,
implicitly, the values of $\alpha_j$ considered in $\mathcal{O}(x_0)$ are such that the associated flow $t \mapsto e^{t h_j}(x)$ is well defined on $[0,\alpha_j]$ at the considered point $x$. These values of $\alpha_j$ may not span the entire real line.

\noindent \emph{Step 1: We prove that $x_0$ is an interior point of $\mathcal{O}(x_0)$.} By Proposition \ref{p:access}, $\mathcal{O}(x_0)$ has a non empty interior: there exist $x_1 \in \mathcal{O}(x_0)$ and $r>0$ such that $B_{\cM}(x_1,r)  \subset \mathcal{O}(x_0)$. By definition of $\mathcal{O}(x_0)$, there exist $n \in \N$, $h_1,\dots, h_n \in F$, $\alpha_1,\dots,\alpha_n \in \R$ such that $x_1=e^{\alpha_n h_n} \dots e^{\alpha_1 h_1}(x_0)$. Since 
$e^{-\alpha_1 h_1} \dots e^{-\alpha_n h_n}$ is a diffeomorphism of $\cM$, then
$\Omega:=e^{-\alpha_n h_n} \dots e^{-\alpha_1 h_1} B_{\cM}(x_1,r) $ is an open neighborhood of $x_0$ in $\mathcal{O}(x_0)$. 

\medskip

\noindent \emph{Step 2: We prove that $\mathcal{O}(x_0)$ is open.} Let $y \in \mathcal{O}(x_0)$. Then $\mathcal{O}(y) \subset \mathcal{O}(x_0)$ thus $\mathcal{O}(x_0)$ contains an open neighborhood of $y$, by Step 1. 

\medskip

\noindent \emph{Step 3: Conclusion.} For $x_0, x_* \in \cM$, $x_* \in \mathcal{O}(x_0)$ is equivalent to $x_0 \in \mathcal{O}(x_*)$. Thus the relation $x_0 \sim x_*$ iff $x_* \in \mathcal{O}(x_0)$ defines an equivalence relation. $\cM$ is the disjoint union of the equivalence classes. By Step 2, these classes are open. Since $\cM$ is connected, $\cM$ is the unique class.
\end{proof}

\subsection{Strong LARC implies STC} \label{subsec:strong_LARC}

\begin{proof}[Proof of Theorem \ref{thm:Old_SC}] 
\noindent \emph{Step 1: We prove the system is small-time approximately controllable.} Let $x_0,x_* \in \cM$ and $\epsilon>0$. By Theorem \ref{thm:Chow}, there exist $n \in \N$, $h_1, \dots,h_n \in \{f_1,\dots, f_m\}$, $\alpha_1,\dots,\alpha_n \in \R$ such that $x_*=e^{\alpha_n h_n} \dots e^{\alpha_1 h_1}(x_0)$. For $t>0$ small enough, the point $e^{tf_0+\alpha_n h_n} \dots e^{tf_0+\alpha_1 h_1}(x_0)$ is accessible from $x_0$ in time $nt$, and converges towards $x_*$ as $t \to 0$.

\medskip

\noindent \emph{Step 2: We prove the system is small-time exactly controllable.}  Let $x_0, x_* \in \cM$ and $\epsilon>0$. We apply Proposition \ref{Prop:Access_0} to the system running backward in time $\dot{x}=-f_0(x)-\sum_{j=1}^m u_j(t) f_j(x)$: there exists an open subset $\Omega_{\epsilon}$ of $\cM$ such that, for every $x_1 \in \Omega_{\epsilon}$, there exist $T \in [0,\epsilon/2]$ and $u \in L^1(0,T)$ such that $x(T;u,x_1)=x_*$. By Step 1, there exists $T' \in [0,\epsilon/2]$, $u' \in L^1(0,T')$ such that $x_1:=x(T';u',x_0)$ belongs to $\Omega_{\epsilon}$. We get the conclusion with the time $T+T' \in [0,\epsilon]$ and the concatenated control $u \diamond u'$.
\end{proof}

\subsection{Estimates on flows} \label{subsec:estim_flow}

We recall that $\cM$ is a smooth submanifold of $\R^d$ and that a smooth vector field on $\cM$ is the restriction to $\cM$ of a smooth vector field on an open neighborhood of $\cM$ in $\R^d$, which is tangent to $\cM$. Thus, all terms in the statements below are well-defined: the sum is the one over $\R^d$, the differentials correspond to those over $\R^d$, etc.


\begin{lemma} \label{Lem_estim_flots}
If $f$ is a $\mathcal{C}^1$ (resp. $\mathcal{C}^2$) vector field on $\cM$ then the estimate
$$e^{tf}(x)=x+tf(x)+\underset{t \to 0}{O}(t^2)
\qquad
\text{(resp. } \quad
e^{tf}(x)=x+tf(x)+\frac{t^2}{2} df(x)f(x) + 
\underset{t \to 0}{O}(t^3) \quad \text{)}
$$
holds uniformly with respect to $x$ in any compact subset of $\cM$.
\end{lemma}

\begin{proof} Let $K$ be a compact subset of $\cM$ and $K'=\left( K+\overline{B}(0,1) \right) \cap \cM$.
By continuity, there exists $t^*>0$ such that $e^{[0,t^*]f}(K) \subset K'$.
For $x \in K$ and $t \in [0,t^*]$, by the Taylor-Lagrange inequality 
$$\left|  e^{tf}(x) - x - t f(x) \right|
\leq \frac{t^2}{2} \underset{\tau \in [0,t]}{\sup} \left| 
df(e^{\tau f}(x))\, f(e^{\tau f}(x)) \right| \leq M_1 t^2,$$
$$\left|  e^{tf}(x) - x - t f(x) - \frac{t^2}{2} df(x) f(x) \right|
\leq  M_2 t^3,$$
where $M_1:=\frac{1}{2} \| df.f \|_{L^{\infty}(K')}$ and
$M_2=\frac{1}{3!} \| d^2 f(f,f) + (df)^2 f\|_{L^{\infty}(K')}$.
\end{proof}

\begin{lemma} \label{Lem_iteration}
Let 
$f$ be a $\mathcal{C}^1$ vector field on $\cM$,
$\sigma \in (0,1]$ and
$(\Phi_t)_{t \in [0,1]}$ be a family of maps 
$\Phi_t:\text{Dom}(\Phi_t) \subset \cM \to \cM$ such that 
for every compact subset $K'$ of $\cM$,
there exists $t_{K'}>0$ such that 
$K' \subset \cap_{t \in [0,t_{K'}]} \text{Dom}(\Phi_t)$ and
the estimate 
$$\Phi_t(x)=x+tf(x)+\underset{t \to 0}{O}(t^{1+\sigma})$$ 
holds uniformly with respect to $x \in K'$.
Then the estimate
$$ \left( \Phi_{\frac{1}{n}} \right)^n(x)=e^{f}(x)+\underset{n \to \infty}{O}(n^{-\sigma})$$ 
holds uniformly with respect to $x$ in any compact subset of $\text{Dom}(e^f)$.
\end{lemma} 

\begin{proof}
Let $K$ be a compact subset of $\text{Dom}(e^f)$.
Then $K' := \left( e^{[0,1]f}(K)+\overline{B}(0,1) \right) \cap \cM$ is a compact subset of $\cM$.
There exists $t_0>0$ such that $K' \subset \text{Dom}(e^{t_0 f})$.
By Gronwall Lemma, 
$$\forall t \in [0,t_0], x,x' \in K', \quad
\left| e^{tf}(x)-e^{tf}(x') \right|
\leq |x-x'| e^{C_0 t} \leq \left( 1+ Ct \right)|x-x'|
$$
where 
$C_0:=\| \nabla f \|_{L^{\infty}(K'')}$,
$K'':=\overline{\text{Conv}}( e^{[0,t_0]}(K') ) \cap \cM$ is a compact subset of $\cM$ and $C:=C_0 e^{C_0 t_0}$. 
By assumption and Lemma \ref{Lem_estim_flots}, there exist constants $t^*,M>0$ such that 
$$\forall t \in [0,t^*], x \in K'
\qquad
\left| \Phi_t(x)-e^{tf}(x) \right| \leq M t^{1+\sigma}.$$
Let $n \in \N^*$ large enough so that 
$n^{-1} < \min\{t_0 ; t^*\}$ and $M n^{-\sigma} e^{C}<1$. 
Let $x \in K$.
The sequence $y_k:=e^{\frac{k}{n}f}(x)$ is well defined for every $k \in \llbracket 0 , n \rrbracket$ because $x \in \text{Dom}(e^f)$; moreover
$y_k \in e^{[0,1]f}(K)$.
We define $z_0=x$ and, by induction, if $z_k \in K'$ then 
$z_{k+1}:=(\Phi_{\frac 1n})(z_k)$.

If $k \in \llbracket 0 , n-1 \rrbracket$ and $z_0,\dots,z_k \in K'$ then, for every $j \in \llbracket 0, k \rrbracket$
$$\left|z_{j+1}-y_{j+1}\right|
\leq
\left| \Phi_{\frac{1}{n}}(z_j)-e^{\frac{1}{n}f}(z_j) \right| + 
\left| e^{\frac{1}{n}f}(z_j) - e^{\frac{1}{n}f}(y_j)\right|
\leq \frac{M}{n^{1+\sigma}} + \left( 1+\frac{C}{n} \right)|y_j-z_j|
$$
i.e. the sequence $\epsilon_j:=\left( 1+\frac{C}{n} \right)^{-j} |y_j-z_j|$ satisfies
$\epsilon_{j+1} \leq \epsilon_j+\frac{M}{n^{1+\sigma}}$ thus
$\epsilon_{k+1} \leq \frac{(k+1) M}{n^{1+\sigma}} \leq M n^{-\sigma}$ and
$$\left|y_{k+1}-z_{k+1}\right| \leq M n^{-\sigma} \left( 1+\frac{C}{n} \right)^k \leq M n^{-\sigma} e^C < 1.$$
This proves that $z_0,\dots,z_n \in K'$ and 
$| z_n-y_n | 
\leq  M n^{-\sigma} e^{C}$ which is the conclusion.
\end{proof}

\subsection{A conjugation formula}  \label{subsec:conj}

\begin{definition} \label{Def:push}
    Let $\theta$ be a smooth diffeomorphism of $\cM$ and $h$ be a smooth vector field on $\cM$. The push-forward of $h$ by $\theta$ is the smooth vector field $(\theta_{\star} h)$ on $\cM$ defined, for any $y \in \cM$ by
    $$\big( \theta_{\star} h \big)(y) := d\theta\big( \theta^{-1}(y) \big) h \big( \theta^{-1}(y) \big).$$
As a consequence, $\theta e^{th} \theta^{-1} = e^{t (\theta_{\star} h)}$ for every $t \in \R$ for which the flows are defined.
\end{definition}

\begin{proof}
Let $x_0 \in \cM$,
$I$ denote the maximal interval of existence of $t \mapsto e^{th}(x_0)$ and, for every $t \in I$, $y(t):= \theta( e^{th}(x_0) )$. 
By the chain rule, for every $t \in I$,
$\dot{y}(t)= d\theta( e^{th}(x_0) ) h (e^{th}(x_0))=(\theta_{\star} h)(y(t))$ thus $y(t)=e^{t(\theta_{\star}h)}y(0)$,
or equivalently $\theta (e^{th}(x_0))=e^{t(\theta_{\star}h)} \theta(x_0)$.
\end{proof}

\begin{lemma} \label{Lem:Conj}
Let $g,h$ be smooth vector fields on $\cM$ such that the vector space $E:=\text{Span}\{ \ad_g^k(h) ; k \in \N \}$ has finite dimension. Then, for every $t>0$ and $x \in \text{Dom}(e^{tg})$ 
   $$ \Big((e^{-tg})_{\star} h\Big)(x)=\sum_{k=0}^{\infty} \frac{t^k}{k!} \ad_g^k(h)(x). $$ 
In particular $e^{-tg} e^{h} e^{tg} = e^{F}$ where $F=\exp(t \ad_g)(h)$.
\end{lemma}

\begin{proof}
Let $t_*>0$, $x \in \text{Dom}(e^{t_*g})$ and $\phi:[0,t_*]\to \R^d$ defined by $\phi(t):= \big((e^{-tg})_{\star} h\big)(x)$. Using the chain rule, we obtain $\phi'(t)=\big((e^{-tg})_{\star} [g,h]\big)(x)$. We deduce by induction that for every $k \in \N$, $\phi^{(k)}(t)=\big((e^{-tg})_{\star} \ad_{g}^k(h)\big)(x)$. In particular $\phi^{(k)}(0)=\ad_{g}^k(h)(x)$.

Let $t \in [0,t_*]$. For $n \in \N^*$, using the Taylor formula, we obtain
\begin{equation} \label{Taylor}
\left| \phi(t)-\sum_{k=0}^{n-1} \frac{t^k}{k} \phi^{(k)}(0) \right|
\leq \frac{|t|^n}{n!} \underset{\tau \in [0,t]}{\sup} 
\left|
\big((e^{-\tau g})_{\star} \ad_{g}^n(h)\big)(x)\right|.
\end{equation}
$K:=e^{[0,t_*]g}(x)$ is a compact subset of $\cM$.
The linear map $\ad_g$ is continuous on the finite-dimensional space $E$ thus there exists $C>0$ such that
$$\forall f \in E, \qquad \| \ad_g(f) \|_{L^{\infty}(K)} \leq C \| f \|_{L^{\infty}(K)} $$
Using this estimate in \eqref{Taylor} gives $C'=C'(x,t_*,g,h)$ such that, for every $t \in [0,t_*]$ and $n \in \N^*$
$$ \left| \Big((e^{-tg})_{\star} h\Big)(x) - \sum_{k=0}^{n-1} \frac{t^k}{k!} \ad_g^k(h)(x) \right|
\leq C' \frac{|Ct|^n}{n!} $$
which proves the convergence of the series.  
\end{proof}
\medskip
\textbf{Acknowledgments.}
The authors are supported by 
ANR-25-CE40-4062 (project QUEST),
ANR-24-CE40-3008-01 (project QuBiCCS),
ANR-24-CE40-5470 (project CHAT), 
ANR-11-LABX-0020 (Centre Henri Lebesgue)
and from the Fondation Simone et Cino Del Duca – Institut de France.

\bibliographystyle{abbrv}
\bibliography{biblio}

\end{document}